\magnification=\magstep1
\input amstex
\documentstyle{amsppt}
\catcode`\@=11 \loadmathfont{rsfs}
\def\mycal{\mathfont@\rsfs}
\csname rsfs \endcsname \catcode`\@=\active

\vsize=6.5in

\topmatter
\title ON THE SUPERRIGIDITY OF MALLEABLE ACTIONS \\
WITH SPECTRAL GAP
\endtitle
\author SORIN POPA \endauthor

\rightheadtext{superrigidity of malleable actions}

\affil University of California, Los Angeles\endaffil

\address Math.Dept., UCLA, LA, CA 90095-155505\endaddress
\email popa\@math.ucla.edu\endemail

\thanks Research supported in part by NSF Grant 0601082.\endthanks

\abstract We prove that if a countable group $\Gamma$ contains
infinite commuting subgroups $H, H'\subset \Gamma$ with $H$
non-amenable and $H'$ ``weakly normal'' in $\Gamma$, then any
measure preserving $\Gamma$-action on a probability space which
satisfies certain malleability, spectral gap and weak mixing
conditions (e.g. a Bernoulli $\Gamma$-action) is cocycle superrigid.
If in addition $H'$ can be taken non-virtually abelian and $\Gamma
\curvearrowright X$ is an arbitrary free ergodic action while
$\Lambda \curvearrowright Y=\Bbb T^\Lambda$ is a Bernoulli action of
an arbitrary infinite conjugacy class group, then any isomorphism of
the associated II$_1$ factors $L^\infty X \rtimes \Gamma \simeq
L^\infty Y \rtimes \Lambda$ comes from a conjugacy of the actions.



\endabstract

\endtopmatter

\document

\heading 1. Introduction \endheading

Some of the most interesting aspects of the dynamics of measure
preserving actions of countable groups on probability spaces,
$\Gamma \curvearrowright (X,\mu)$, are revealed by the study of {\it
group measure space} von Neumann algebras $L^\infty(X)\rtimes
\Gamma$ ([MvN1]) and the classification of actions up to {\it orbit
equivalence} ({\it OE}), i.e. up to isomorphism of probability
spaces carrying the orbits of actions onto each other. Although one
is in von Neumann algebras and the other in ergodic theory, the two
problems are closely related, as an OE of actions $\Gamma
\curvearrowright X$, $\Lambda \curvearrowright Y$ has been shown to
implement an algebra isomorphism $L^\infty(X)\rtimes \Gamma \simeq
L^\infty(Y)\rtimes \Lambda$ taking $L^\infty(X)$ onto $L^\infty(Y)$,
and vice-versa ([Si], [Dy], [FM]). In particular, the isomorphism
class of $L^\infty(X)\rtimes \Gamma$ only depends on the equivalence
relation $\Cal R_\Gamma=\{(t,gt) \mid t\in X, g\in \Gamma\}$.

Thus, {\it von Neumann equivalence} ({\it vNE}) of group actions,
requiring isomorphism of their group measure space algebras, is
weaker than OE. Since there are examples of
non-OE actions whose associated von Neumann algebras
are all isomorphic ([CJ1]), it is in general strictly weaker. On the
other hand, OE is manifestly weaker than classical {\it conjugacy},
which for free actions $\Gamma \curvearrowright X$, $\Lambda
\curvearrowright Y$ requires isomorphism of probability spaces
$\Delta: (X,\mu) \simeq (Y,\nu)$ satisfying $\Delta \Gamma
\Delta^{-1}=\Lambda$ (so in particular $\Gamma\simeq \Lambda$). How
much weaker vNE and OE can be with respect to conjugacy is best seen in the
amenable case, where by a celebrated theorem of Connes all free
ergodic actions of all (infinite) amenable groups give rise to the
same II$_1$ factor ([C1]) and by ([Dy], [OW], [CFW]) they are
undistinguishable under OE as well. Also, any embedding of algebras
$L^\infty(X)\rtimes \Gamma \subset L^\infty(Y)\rtimes \Lambda$ with
$\Lambda$ amenable forces $\Gamma$ to be amenable.

But the non-amenable case is extremely complex, and for many years
progress has been slow ([MvN2], [Dy], [Mc], [C2], [CW], [Sc]), even
after the discovery of the first rigidity phenomena, by Connes in
von Neumann algebras ([C3,4]) and by Zimmer in OE ergodic theory
([Z1,2]). This changed dramatically over the last 7-8 years, with
the advent of a variety of striking rigidity results ([Fu1], [G1,2],
[MoSh], [P1-8], [H], [HK], [Ki]; see [P9] for a survey; also [Sh] for a
survey on OE rigidity).

Our aim in this paper is to investigate the most ``extreme'' such
phenomena, called {\it strong rigidity}, which show that for certain
classes of {\it source} group actions $\Gamma \curvearrowright X$
and {\it target} actions $\Lambda \curvearrowright Y$ any
isomorphism $L^\infty(X)\rtimes \Gamma \simeq L^\infty(Y)\rtimes
\Lambda$ (resp. any OE of $\Gamma \curvearrowright X$, $\Lambda
\curvearrowright Y$) comes from a conjugacy, modulo perturbation by
an inner automorphism of $L^\infty(Y)\rtimes \Lambda$ (resp. of
$\Cal R_\Lambda$). Ideally, one seeks to prove this under certain
conditions on the source group actions $\Gamma \curvearrowright X$
but no condition at all (or very little) on the target $\Lambda
\curvearrowright Y$, a type of result labeled {\it superrigidity}.
On the orbit equivalence side, such results appeared first in [Fu1]
(for actions of higher rank lattices, such as $SL(n,\Bbb Z)
\curvearrowright \Bbb T^n$, $n \geq 3$) and then in [MoSh] (for
doubly ergodic actions of products of word hyperbolic groups, such
as $\Bbb F_n \times \Bbb F_m$). In the meantime, new developments in
von Neumann algebras ([P3,8]) led to the first vNE strong rigidity
result in [P4,5]. It shows that any isomorphism of group measure
space factors $L^\infty(X)\rtimes \Gamma \simeq L^\infty(Y)\rtimes
\Lambda$, with $\Gamma$ an {\it infinite conjugacy class} ({\it
ICC}) group having an infinite subgroup satisfying a weak normality
condition and with the relative property (T) of Kazhdan-Margulis
($\Gamma$ {\it w-rigid}) and $\Lambda \curvearrowright Y$ a
Bernoulli action of an arbitrary ICC group, comes from a conjugacy.
While obtained in a purely von Neumann algebra framework, this
result provides new OE rigidity phenomena as well, showing for
instance that Bernoulli actions of Kazhdan groups are OE superrigid
([P5]).

The ideas and techniques in [P3,4,5] were further exploited in
[P1] to obtain a cocycle superrigidity result for Bernoulli
actions of w-rigid groups $\Gamma$, from which OE superrigidity is
just a consequence. Thus, [P1] shows that any measurable cocycle
for $\Gamma \curvearrowright X=[0,1]^\Gamma$ is $\Cal U_{fin}$-{\it
cocycle superrigid} ({\it CSR}), i.e. any $\Cal V$-valued cocycle
for $\Gamma \curvearrowright X$ is cohomologous to a group morphism
$\Gamma \rightarrow \Cal V$, whenever $\Cal V$ is a closed subgroup
of the unitary group of a separable finite von Neumann algebra, for
instance if $\Cal V$ is countable discrete, or separable compact.

The sharp OE and vNE rigidity results in [P1,5], and in fact in
[P1-8], [PS], [IPeP], [PV], [V], [I1,2] as well, are due to a
combination (co-existence) of deformability and rigidity assumptions
on the group actions. The deformability condition imposed is often
the {\it malleability} of the action (e.g. in [P1,5]), a typical
example of which are the Bernoulli actions, while the rigidity
assumption is each time some weak form of property (T) (on the
acting group, as in [P1,5], or on the way it acts, as in [P8]).
Thus, the {\it deformation/rigidity} arguments used in all these
papers seemed to depend crucially on the ``property (T)-type''
assumption.

However, in this paper we succeed to remove this assumption
completely. Namely, we prove a new set of rigidity results for
malleable actions, in some sense ``parallel'' to the ones in
[P1,5], but which no longer assume Kazhdan-type conditions on the
source group, being surprisingly general in this respect. For
instance, we show that if $\Cal V\in \Cal U_{fin}$ and $\Gamma$ is
an arbitrary group, then any $\Cal V$-valued cocycle for a Bernoulli
$\Gamma$-action can be untwisted on the centralizer (or commutant)
of any non-amenable subgroup $H$ of $\Gamma$! More precisely, we
prove (compare with 5.2/5.3 in [P1]):

\proclaim{1.1. Theorem (CSR: s-malleable actions)} Let
$\Gamma \curvearrowright (X,\mu)$ be an m.p. action
of a countable group $\Gamma$. Let $H,H'\subset \Gamma$
be infinite commuting subgroups
such that:

\vskip .05in

$(a)$ $H \curvearrowright X$ has stable spectral gap.

$(b)$ $H' \curvearrowright X$ is weak mixing.

$(c)$ $HH' \curvearrowright X$ is s-malleable.

\vskip .05in Then $\Gamma \curvearrowright X$ is $\mycal
U_{fin}$-cocycle superrigid on $HH'$. If in addition $H'$ is
w-normal in $\Gamma$, or $H'$ is wq-normal but $\Gamma
\curvearrowright X$ is mixing, then $\Gamma \curvearrowright X$ is
$\mycal U_{fin}$-cocycle superrigid on all $\Gamma$. Moreover, the
same conclusions hold true if we merely assume $HH'\curvearrowright
X$ to be a relative weak mixing quotient of an  m.p. action $HH'
\curvearrowright (X',\mu')$ satisfying conditions $(a),(b),(c)$.
\endproclaim

The {\it stable spectral gap} condition $(a)$
in Theorem $1.1$ means
the representation implemented by the action $H \curvearrowright X$
on $L^2 X \overline{\otimes} L^2 X \ominus \Bbb C$ has spectral gap,
i.e. has no approximately invariant vectors (see Section 3). It
automatically  implies $H$ is non-amenable. The {\it s-malleability}
condition for an m.p. action $\Gamma_0 \curvearrowright X$ was
already considered in [P1-5] and is discussed in Section 2. An
action $\Gamma_0 \curvearrowright (X,\mu)$ is a {\it relative weak
mixing} quotient of an m.p. action $\Gamma_0 \curvearrowright
(X',\mu')$ if it is a quotient of it and $\Gamma_0 \curvearrowright
X'$ is weak mixing relative to $\Gamma_0 \curvearrowright X$ in the
sense of [F], [Z3] (see also Definition 2.9 in [P1]).

The two ``weak normality'' conditions considered in
Theorem 1.1 are the same as
in [P1,2,5]: An infinite subgroup $\Gamma_0\subset \Gamma$ is {\it
w-normal} (resp. {\it wq-normal}) in $\Gamma$ if there exists a well
ordered family of intermediate subgroups $\Gamma_0 \subset \Gamma_1
\subset ... \subset \Gamma_{\jmath} \subset ... \subset
\Gamma_{\imath} = \Gamma$ such that for each $0 < \jmath \leq
\imath$, the group $\Gamma'_\jmath=\bigcup_{n < \jmath} \Gamma_\jmath$
is normal in $\Gamma_\jmath$ (resp. $\Gamma'_{\jmath}$ is generated
by the elements $g \in \Gamma$ with $|g\Gamma'_{\jmath}g^{-1} \cap
\Gamma'_{\jmath}|=\infty$).

Any generalized Bernoulli action $\Gamma_0 \curvearrowright \Bbb
T^I$, associated to an action of a countable group $\Gamma_0$ on a
countable set $I$, is s-malleable. Given any probability space
$(X_0,\mu_0)$ (possibly atomic), the generalized Bernoulli action
$\Gamma_0 \curvearrowright (X_0,\mu_0)^I$ is a relative mixing
quotient of the s-malleable action $\Gamma_0 \curvearrowright \Bbb
T^I$. Any Gaussian action $\sigma_\pi:\Gamma_0 \curvearrowright
(\Bbb R, (2\pi)^{-1/2} \int \cdot e^{-t^2}\text{\rm d} t)^n$
associated to an orthogonal representation $\pi$ of $\Gamma_0$ on
the $n$-dimensional real Hilbert space $\Cal H_n=\Bbb R^n$, $2\leq n
\leq \infty$, is easily seen to be s-malleable (cf.
[Fu2]). The action $\sigma_\pi$ has stable
spectral gap on some subgroup $H\subset \Gamma_0$ once the
orthogonal representation $\pi_{|H}$ has stable spectral gap. By
[P2], a sufficient condition for a generalized Bernoulli action $H
\curvearrowright (X_0,\mu_0)^I$ to have stable spectral gap is that
$\{g \in H \mid gi=i\}$ be amenable, $\forall i\in I$. Thus, Theorem
1.1 implies:

\proclaim{1.2. Corollary (CSR: Bernoulli and Gaussian actions)} Let
$\Gamma$ be a countable group having infinite commuting subgroups
$H,H'$ with $H$ non-amenable. Let $\Gamma \curvearrowright X$ be an
m.p. action whose restriction to $HH'$ is a relative weak mixing
quotient of one of the following: \vskip .05in $1^\circ.$ A
generalized Bernoulli action $HH' \curvearrowright (X_0,\mu_0)^I$,
with the actions of $H,H'$ on the countable set $I$ satisfying $|H'
i|=\infty$ and $\{g \in H \mid gi=i\}$ amenable, $\forall i\in I$.

$2^\circ.$ A Gaussian action associated to an orthogonal
representation of $HH'$ which has stable spectral gap on $H$ and no
finite dimensional $H'$-invariant subspaces. \vskip .05in

If $H'$ is w-normal in $\Gamma$, then
$\Gamma \curvearrowright X$ is $\mycal U_{fin}$-cocycle superrigid.
If $H'$ is merely wq-normal in $\Gamma$ but $\Gamma \curvearrowright
X$ is a weak mixing quotient of a Bernoulli action, then again
$\Gamma \curvearrowright X$ is $\mycal U_{fin}$-cocycle superrigid.
\endproclaim

Due to Theorems 5.6-5.8 in [P1], the cocycle superrigidity
results in Theorem 1.1 and Corollary 1.2
imply several superrigidity results in orbit equivalence ergodic
theory:

\proclaim{1.3. Corollary (OE superrigidity results)} Let $\Gamma$ be
a countable group with no finite normal subgroups and having
infinite commuting subgroups $H,H'$, with $H$ non-amenable. Assume
the free m.p. action $\Gamma \curvearrowright X$ is a relative weak
mixing quotient of an s-malleable action $\Gamma \curvearrowright
(X',\mu')$ such that:

\vskip .05in \noindent $(1.3)$ $H
\curvearrowright X'$ has stable spectral gap and either $H'$ is
w-normal in $\Gamma$ with $H' \curvearrowright X'$ weak mixing, or
$H$ is merely wq-normal in $\Gamma$ but with $\Gamma
\curvearrowright X'$ mixing.

\vskip .05in Then $\Gamma  \curvearrowright X$ satisfies the
conclusions in $5.6, 5.7, 5.8$ of $\text{\rm [P1]}$. In
particular, any Bernoulli action $\Gamma \curvearrowright
(X,\mu)=(X_0,\mu_0)^\Gamma$ is OE superrigid, i.e. any OE between
$\Gamma \curvearrowright X$ and an arbitrary free ergodic m.p.
action of a countable group comes from a conjugacy. If in addition
$\Gamma$ is ICC, then ANY quotient of $\Gamma \curvearrowright X$
which is still free on $\Gamma$ is OE superrigid.
\endproclaim

By 2.7 in [P2], Corollary 1.2 provides a large class of groups
with uncountably many OE inequivalent actions, adding to the
numerous examples already found in [Z2], [H], [GP], [P2], [I1],
[IPeP]:

\proclaim{1.4. Corollary} Assume $\Gamma$ contains a non-amenable
subgroup with its centralizer infinite and wq-normal in $\Gamma$, e.g.
$\Gamma$ non-amenable and either a product of two infinite groups,
or having infinite center. Then for any countable abelian group $L$
there exists a free ergodic action $\sigma_L$ of $\Gamma$ on a
probability space with the first cohomology group $\text{\rm
H}^1(\sigma_L)$ equal to $\text{\rm Hom}(\Gamma, \Bbb T) \times L$. In
particular, $\Gamma$ has a continuous family of OE inequivalent
actions.
\endproclaim

The cocycle superrigidity result
in Theorem 1.1 is analogous to 5.2/5.3 in [P1]. The
trade-off for only assuming $H\subset \Gamma$ non-amenable in
Theorem 1.1,
rather than Kazhdan (as in [P1]), is the spectral gap condition on
the action. The proof is still based on a deformation/rigidity
argument, but while the malleability is combined in [P1] with
property (T) rigidity, here it is combined with spectral gap
rigidity. Also, rather than untwisting a given cocycle on $H$, we
first untwist it on the group $H'$ commuting with $H$. Due to the
weak mixing property $(b)$ in Theorem $1.1$ and 3.6 in [P1], it then gets
untwisted on the w-normalizer of $H'$, thus on $HH'$. Altogether, we
rely heavily on technical results from [P1].

\vskip .05in

We use the same idea of proof, combined this time with technical
results from [P4,5], to obtain a vNE strong rigidity result
analogue to 7.1/7.1' in [P5], which derives conjugacy of actions
from the isomorphism of their group measure space factors. Note that
while the ``source'' group $\Gamma$ is still required to have a
non-amenable subgroup with infinite centralizer, the $\Gamma$-action
here is completely arbitrary. In turn, while the ``target'' group
$\Lambda$ is arbitrary, the $\Lambda$-action has to be Bernoulli.
Thus, the spectral gap condition, which is automatic for Bernoulli
actions, is now on the target side.

\proclaim{1.5. Theorem (vNE strong rigidity)} Assume $\Gamma$
contains a non-amenable subgroup with centralizer non-virtually
abelian and wq-normal in $\Gamma$. Let $\Gamma \curvearrowright
(X,\mu)$ be an arbitrary free ergodic m.p. action. Let $\Lambda$ be
an arbitrary ICC group and $\Lambda \curvearrowright (Y,\nu)$ a
free, relative weak mixing quotient of a Bernoulli action $\Lambda
\curvearrowright (Y_0,\nu_0)^\Lambda$. If $\theta: L^\infty X
\rtimes \Gamma \simeq (L^\infty Y \rtimes \Lambda)^t$ is an
isomorphism of $\text{\rm II}_1$ factors, for some $0< t \leq 1$,
then $t=1$ and $\theta$ is of the form $\theta = {\text{\rm Ad}}(u)
\circ \theta^\gamma \circ \theta_0$, where: $u$ is a unitary element
in $L^\infty Y \rtimes \Lambda$; $\theta^\gamma \in {\text{\rm
Aut}}(L^\infty Y \rtimes \Lambda)$ is implemented by a character
$\gamma$ of $\Lambda$; $\theta_0 :L^\infty X \rtimes \Gamma \simeq
L^\infty Y \rtimes \Lambda$ is implemented by a conjugacy of $\Gamma
\curvearrowright X$, $\Lambda \curvearrowright Y$.
\endproclaim

When applied to isomorphisms $\theta$ implemented by OE of actions,
Theorem 1.5 above implies an {\it OE Strong Rigidity} result analogue to
7.6 in [P5]. We in fact derive an even stronger rigidity result,
for embeddings of equivalence relations, parallel to 7.8 in [P5]:

\proclaim{1.6. Theorem (OE strong rigidity for embeddings)} Let
$\Gamma \curvearrowright (X,\mu)$, $\Lambda \curvearrowright
(Y,\nu)$  be as in $1.5$. If $\Delta:(X,\mu)\simeq (Y,\nu)$ takes
each $\Gamma$-orbit into a $\Lambda$-orbit (a.e.), then there exists
a subgroup $\Lambda_0\subset \Lambda$ and $\alpha\in \text{\rm
Inn}(\Cal R_\Lambda)$ such that $\alpha\circ \Delta$ conjugates
$\Gamma \curvearrowright X$, $\Lambda_0 \curvearrowright Y$.
\endproclaim

Notice that
although OE superrigidity  results are of a stronger type than OE
strong rigidity, Theorem 1.6 cannot be deduced from
Corollary 1.3, nor in fact
from the cocycle superrigidity result in Theorem 1.1. Likewise,
the OE strong rigidity
7.8 in [P5] cannot be derived from results in [P1].

The idea of combining malleability with spectral gap rigidity in the
proofs of Theorems 1.1 and 1.5 is inspired from [P6], where a
similar argument was used to prove a unique ``McDuff decomposition''
result for II$_1$ factors. We comment on this in Section 6, were we
include other remarks as well, showing for instance that if $\Gamma
\curvearrowright X$ is a Bernoulli action of a non-amenable group,
then $L^\infty X \rtimes \Gamma$ follows prime, due to the same
arguments. We also revisit the Connes-Jones counterexample in
[CJ1] and point out that, due to results in [P1,5] and in this
paper, it provides cocycle superrigid (in particular OE superrigid)
actions $\Gamma \curvearrowright X$ whose equivalence relation $\Cal
R_\Gamma$ has trivial fundamental group, $\Cal F(\Cal
R_\Gamma)=\{1\}$, while the associated II$_1$ factor
$M=L^\infty(X)\rtimes \Gamma$ has fundamental group equal to $\Bbb
R_+^*$, so $M$ can be realized by uncountably many OE-inequivalent
actions (two of which are free).
Section 4 contains the proof of
Theorem 1.1 and Corollaries 1.2, 1.3, 1.4 and Section 5 the
proof of Theorems
1.5 and 1.6. Both sets of proofs rely heavily on technical
results from [P1] and respectively [P4,5]. The present paper
should in fact be viewed as a companion to [P1,4,5], from which
notations and terminology are taken as well.

In Section 2 we comment on s-malleability of actions and
transversality, then in Section 3 we define the notion of stable
spectral gap for actions and representations of groups, and examine
how Bernoulli, Gaussian and Bogoliubov actions (which are the basic
examples of s-malleable actions) behave with respect to this
property.

I am extremely grateful to Stefaan Vaes and the referee for pointing
out to me a redundancy in Section 2 of the initial version of the
paper (see Section 6.2).

\heading 2. Transversality of s-malleable actions \endheading

In [P1-4] we have considered various degrees of malleability for
actions of groups on probability measure spaces $\Gamma
\curvearrowright (X,\mu)$ (more generally on von Neumann algebras).
The weakest such condition (2.1 in [P2], or 4.2 in [P1]) requires
the connected component of the identity in the centralizer
${\text{\rm Aut}}_\Gamma (X\times X, \mu \times \mu)$ of the double
action $g(t,t')=(gt,gt')$, $\forall (t,t')\in X \times X$, $g\in
\Gamma$, to contain an automorphism $\alpha_1$ satisfying
$\alpha_1(L^\infty X \otimes 1)=1 \otimes L^\infty X$, when viewed
as an automorphism of function algebras. More generally, an action
on a finite von Neumann algebra $\Gamma \curvearrowright (P,\tau)$
is malleable if it admits an extension to an action on a larger
finite von Neumann algebra, $\Gamma \curvearrowright
(\tilde{P},\tilde{\tau})$, such that the  connected component of $id$
in the centralizer of this action, ${\text{\rm Aut}}_\Gamma
(\tilde{P}, \tilde{\tau})$, contains an automorphism $\alpha_1$ with
$P_1=\alpha_1(P)$ perpendicular to $P$ (with respect to
$\tilde{\tau}$) and sp$PP_1$ dense in $L^2(\tilde{P})$, in other
words $L^2(\tilde{P})=L^2(\text{\rm sp}PP_1)\simeq
L^2(P)\overline{\otimes} L^2(P_1)$ (see 1.4, 1.5 in [P4]).

It is this condition that we will generically refer to as (basic)
{\it malleability}. We mention that in all existing examples of
malleable actions $\alpha_1$ can in fact be chosen to be the {\it
flip} $(t,t') \mapsto (t',t)$.

A stronger form of malleability in [P1,3,4] requires that there
actually is a continuous group-like ``path'' between the identity
and $\alpha_1$, i.e. a continuous action $\alpha$ of the reals  on
$(X \times X, \mu \times \mu)$, commuting with $\Gamma
\curvearrowright X\times X$ (resp. $\alpha(\Bbb R) \subset
{\text{\rm Aut}}_\Gamma (\tilde{P}, \tilde{\tau})$) and such that
$\alpha_1(L^\infty X \otimes 1)=1\otimes L^\infty X$ (resp.
$\alpha_1(P)=P_1 \perp P$ as before). We call such a path a {\it
malleable deformation} (or {\it path}) of $\Gamma \curvearrowright
X$. An action having such a deformation was still called {\it
malleable} in [P1-4], but to distinguish it from the above weaker
form here we call it {\it path-malleable}.

The strongest condition of this kind considered in [P1-4] requires
the $\pm$ directions on the path $\alpha$ to be symmetric with
respect to the first coordinate of the double space $X \times X$, a
rather natural ``geometric'' property. In rigorous terms, this means
the existence of a period-2 m.p. automorphism $\beta$ of $X \times
X$ commuting with the double $\Gamma$-action  (resp. $\beta
\in {\text{\rm Aut}}_\Gamma (\tilde{P}, \tilde{\tau})$), which acts
as the identity on the first variable (so $\beta (a \otimes
1)=a\otimes 1$, $\forall a\in L^\infty X$; in general $P\subset
\tilde{P}^\beta$) and ``reverses the direction'' of the path
$\alpha$, i.e. $\beta \alpha_t \beta= \alpha_{-t} $, $\forall t$.
Note that $(\alpha,\beta)$ generate a copy of the group of
isometries of $\Bbb R$, Isom$(\Bbb R)=\Bbb R \rtimes \Bbb Z/2\Bbb Z$,
in the centralizer of the double $\Gamma$-action.

This is called {\it s-malleability} in [P1,3,4], and is a useful
strengthening of basic malleability in a ``non-commutative
environment'', e.g. when the probability space is non-commutative
(i.e. $\Gamma$ acts on a finite von Neumann algebra with a trace
$(P,\tau)$), as in [P3], or when malleability is being used to get
information on the von Neumann algebra $L^\infty X \rtimes \Gamma$
and its subalgebras, as in [P4]. Such $(\alpha,\beta)$ plays the
role of a ``device for patching incremental intertwiners'', along
the path $\alpha$. We call the pair $(\alpha,\beta)$ a {\it
s-malleable deformation} (or {\it path}).
One should mention that all known examples
of malleable deformations of actions (generalized Bernoulli actions
[P1-4], Bogoliubov actions [P3] and Gaussian actions [Fu2]) have a
natural symmetry $\beta$ and are thus s-malleable.

Let us note that symmetric deformations automatically satisfy a
natural ``transversality'' condition:

\proclaim{2.1. Lemma} Let $\Gamma \curvearrowright (P,\tau)$ be an
s-malleable action and $(\alpha,\beta)$ the corresponding
s-malleable deformation. Then given any finite von Neumann algebra
$(N,\tau)$ the action $\alpha' = \alpha \otimes 1$ of $\Bbb R$ on
$L^2 \tilde{P} \overline{\otimes} L^2 N$ satisfies \vskip .05in
\noindent $(2.1.1)$. $\|\alpha'_{2s}(x)-x \|_2 \leq 2 \|\alpha'_s(x)-
E_{P\overline{\otimes} N} (\alpha'_s(x)) \|_2$.

\endproclaim
\noindent {\it Proof}. Let $\beta'$ denote the period 2 automorphism
of $L^2 \tilde{P} \overline{\otimes} L^2 N$ given by $\beta'=\beta
\otimes 1$. Since $\beta(x)=x$ for $x\in P$ we have $\beta'(x)=x$
for $x\in L^2P \overline{\otimes} L^2N$. In particular
$\beta'(E_{P\overline{\otimes} N} (\alpha'_s(x)))=
E_{P\overline{\otimes} N} (\alpha'_{s}(x))$. Also, $\beta' \alpha'_t
\beta' = \alpha'_{-t}.$ Thus
$$
\|\alpha'_s(x)-E_{P\overline{\otimes} N} (\alpha'_s(x)) \|_2
=\|\beta'(\alpha'_s(x)-E_{P\overline{\otimes} N} (\alpha'_s(x)))
\|_2
$$
$$
=\|\alpha'_{-s} (x)-E_{P\overline{\otimes} N} (\alpha'_{s}(x))\|_2,
$$
implying that
$$
\| \alpha'_{2s}(x) - x \|_2 = \|\alpha'_s(x)-\alpha'_{-s}(x) \|_2
$$
$$
\leq \|\alpha'_{s}(x)- E_{P\overline{\otimes} N} (\alpha'_s(x)) \|_2
+ \|\alpha'_{-s}(x)- E_{P\overline{\otimes} N} (\alpha'_{s}(x)) \|_2
$$
$$
= 2 \|\alpha'_{s}(x)- E_{P\overline{\otimes} N} (\alpha'_s(x)) \|_2.
$$

\hfill $\square$

\heading 3. Stable spectral gap \endheading

\vskip .05in \noindent {\bf 3.1. Definition}. $1^\circ$.
A unitary representation $\Gamma \curvearrowright \Cal H$ has {\it
spectral gap} (resp. {\it stable spectral gap}) if $1_\Gamma
\nprec \Cal H$ (resp. $1_\Gamma \nprec \Cal H \overline{\otimes}
{\Cal H}^*$). An orthogonal representation has spectral gap
(resp. stable spectral gap) if its complexification has the property.

$2^\circ$. An m.p. action $\Gamma \curvearrowright X$ on a
probability space $(X,\mu)$ has {\it spectral gap} (resp. {\it
stable spectral gap}) if the associated representation $\Gamma
\curvearrowright L^2 X\ominus \Bbb C$ has spectral gap (resp. stable
spectral gap). More generally, if $(P,\tau)$ is a finite von Neumann
algebra, an action $\Gamma \curvearrowright (P,\tau)$ has {\it
spectral gap} (resp. {\it stable spectral gap}) if the
representation $\Gamma \curvearrowright L^2 P\ominus \Bbb C$ has
spectral gap (resp. stable spectral gap).

\proclaim{Lemma 3.2} A representation $\Gamma  \curvearrowright^\pi
\Cal H$ has stable spectral gap if and only if given any
representation $\Gamma \curvearrowright^\rho \Cal K$, the product
representation $(\pi\otimes \rho)(g)=\pi(g)\otimes \rho(g)$ of
$\Gamma$ on $\Cal H\overline{\otimes} \Cal K$ has spectral gap.
\endproclaim
\noindent {\it Proof}. Given Hilbert spaces $\Cal H$, $\Cal K$, we
identify the tensor product Hilbert space $\Cal K \overline{\otimes}
\Cal H^*$ with the Hilbert space $HS(\Cal H, \Cal K)$ of linear
bounded operators $T: \Cal H \rightarrow \Cal K$  of finite
Hilbert-Schmidt norm $\|T\|_{HS}=Tr_{\Cal H}(T^*T)^{1/2}=Tr_{\Cal
K}(TT^*)^{1/2} <\infty$, via the map $\xi \otimes \eta^* \mapsto
T_{\xi \otimes \eta^*}$, with $T_{\xi\otimes\eta^*}(\zeta) =\langle
\zeta, \eta \rangle \xi$, $\zeta \in \Cal H$. In particular, we
identify $\Cal H\overline{\otimes} \Cal H^*$ with the Hilbert space
$HS(\Cal H)$ of Hilbert-Schmidt operators on $\Cal H$.

Note that by the Powers-St\o rmer inequality ([PoSt]),
if for $T,S \in HS(\Cal H, \Cal K)$ we denote $|T|=(T^*T)^{1/2}, |S|
=(S^*S)^{1/2} \in HS(\Cal H)$, then we have the estimate
$\| |T| - |S| \|^2_{HS} \leq \|T-S\|_{HS} \|T+S\|_{HS}$.

Let us now take representations $\Gamma  \curvearrowright^\pi \Cal H$, $\Gamma
\curvearrowright^\rho \Cal K$ and still denote by
$\pi$ the representation of $\Gamma$ on $\Cal H^*$ given by
$\pi(\eta^*)=\pi(\eta)^*$ and by $\rho \otimes \pi$ the
representation of $\Gamma$ on $HS(\Cal H, \Cal K)$ resulting from
the identification of this Hilbert space with $\Cal
K\overline{\otimes} \Cal H^*$. Notice that if $T\in HS(\Cal H, \Cal
K)$, then $(\rho_g \otimes \pi_g(T))^* (\rho_g \otimes \pi_g(T)) =
\pi_g(T^*T)$. From the above Powers-St\o rmer inequality we thus
get

$$
\|\pi_g(|T|) - |T| \|^2_{HS}
$$
$$
\leq\|T+\pi_g \otimes \rho_g (T)\|_{HS} \|T-\pi_g \otimes \rho_g
(T)\|_{HS}
$$
$$
\leq 2 \|T\|_{HS} \|T-\pi_g \otimes \rho_g (T)\|_{HS},
$$
showing that if $T$ is almost invariant to $\rho_g \otimes \pi_g$,
for $g$ in a finite subset $F \subset \Gamma$, then $|T|$ is almost
invariant to $\pi_g, g\in F$. In other words, if $1_\Gamma \prec
\Cal K \overline{\otimes} \Cal H^*$, then $1_\Gamma \prec \Cal H
\overline{\otimes} \Cal H^*$.

\hfill $\square$

\proclaim{3.3. Lemma} $1^\circ$. An orthogonal representation
$\Gamma \curvearrowright^\pi \Cal H_\infty$ on the infinite
dimensional real Hilbert space $\Cal H_\infty$ has stable spectral
gap iff the associated Gaussian (resp. Bogoliubov)
action $\sigma_\pi$ has stable
spectral gap.

$2^\circ$. If $\Gamma$ is non-amenable and $\Gamma \curvearrowright
I$ is so that $\{g\in \Gamma \mid g i = i\}$ is amenable $\forall
i\in I$, then the generalized Bernoulli action $\Gamma
\curvearrowright (X_0,\mu_0)^I$ has stable spectral gap.
\endproclaim
\noindent {\it Proof}. $1^\circ$. For Gaussians, this is clear from
the fact that,  as a representation on $L^2(\Cal
H_\infty,\mu_\infty)$, $\sigma_\pi$ is equivalent to the
representation $\bigoplus_{n \geq 0} \pi_{\Bbb C}^{\odot_{_{s}} n}$,
where $\pi_{\Bbb C}$ is the complexification of $\pi$
and for a representation $\rho$ on a (complex) Hilbert space
$\Cal K$, $\rho^{\odot_{_{s}} n}$ denotes its
$n$'th symmetric tensor power (see e.g.
[CCJJV]). Similarly for Bogoliubov actions.

$2^\circ$. This is Lemma 1.6.4 in [P2]. \hfill $\square$ \vskip
.05in

We mention that in the proof of Theorem 1.1 we will in fact need a weaker
condition on an action $\Gamma \curvearrowright (P,\tau_P)$ than
stable spectral gap, namely a ``stable'' version of the {\it strong
ergodicity} in [Sc], which we recall requires any asymptotically
$\Gamma$-invariant sequence $(x_n) \in (P)_1$ (i.e. $\lim_n \|g x_n
- x_n\|_2 = 0$, $\forall g\in \Gamma$) to be asymptotically scalar
(i.e. $\lim_n \|x_n - \tau(x_n)1\|_2=0$).

\vskip .1in \noindent {\bf 3.4. Definition}. An action $\Gamma
\curvearrowright (P,\tau_P)$ is {\it stably strongly ergodic} if
given any action $\Gamma \curvearrowright (Q,\tau_Q)$ on a finite
von Neumann algebra, any asymptotically $\Gamma$-invariant sequence
of the product action $\Gamma \curvearrowright (P\overline{\otimes}
Q, \tau_P \otimes \tau_Q)$ is (asymptotically) contained in $Q$. An
m.p. action $\Gamma \curvearrowright (X,\mu)$ is {\it stably
strongly ergodic} if the action it implements on $L^\infty X$ is
stably strongly ergodic.

\proclaim{3.5. Lemma} If $\Gamma  \curvearrowright (P,\tau_P)$ has
(stable) spectral gap, then it is (stably) strongly ergodic.
\endproclaim
\noindent
{\it Proof}. This is trivial from the definitions.
\hfill $\square$

\heading 4. Proof of Cocycle Superrigidity \endheading

We use in this section the framework and technical results from
[P1]. Notations that are not specified, can be found in [P1] as
well. We in fact prove a generalized version of Theorem 1.1, for
actions of groups on arbitrary finite von Neumann algebras, which is
the analogue of 5.5 in [P1]. Recall in this respect that if
$\Gamma$ is a discrete group, $\Cal N$ is a finite von Neumann
algebra and $\sigma : \Gamma \rightarrow \text{\rm Aut}(\Cal N)$ an
action of $\Gamma$ on $\Cal N$ (i.e. a group morphism of $\Gamma$
into the group of automorphisms $\text{\rm Aut}(\Cal N)$ of the von
Neumann algebra $\Cal N$), then a (left) {\it cocycle} for $\sigma$
is a map $w: \Gamma \rightarrow \Cal U(\Cal N)$ satisfying $w_g
\sigma_g(w_h)=w_{gh}$, $\forall g,h\in \Gamma$. Also, two such
cocycles $w,w'$ are equivalent if there exists a unitary element
$u\in \Cal U(\Cal N)$ such that $u^* w_g \sigma_g(u)=w'_g$, $\forall
g\in \Gamma$.

\proclaim{4.1. Theorem (Cocycle superrigidity: the non-commutative
case)} Let $\Gamma\curvearrowright^{\sigma_0} (P,\tau)$ be an action
of $\Gamma$ on a finite von Neumann algebra. Let $H,H'\subset
\Gamma$ be infinite commuting subgroups such that:

\vskip .05in

$(a)$ $H \curvearrowright P$ has stable spectral gap.

$(b)$ $H' \curvearrowright P$ is weak mixing.

$(c)$ $HH' \curvearrowright P$ is s-malleable.

\vskip .05in

Let $(N,\tau)$ be an arbitrary finite von Neumann algebra and $\rho$
an action of $\Gamma$ on $(N,\tau)$. Then any cocycle $w$ for the
diagonal product action $\sigma_0\otimes \rho$ of $\Gamma$ on $P
\overline{\otimes} N$ is equivalent to a cocycle $w'$ whose
restriction to $HH'$ takes values in $N=1\otimes N$. If in addition
$H'$ is w-normal in $\Gamma$, or if $H'$ is wq-normal but $\sigma$
is mixing, then $w'$ takes values in $N$ on all $\Gamma$.

Moreover, the same result holds true if $\sigma_0$ extends to an
s-malleable action $\Gamma \curvearrowright^{\sigma_0'} (P',\tau)$
which satisfies $(a), (b), (c)$ and such that $\sigma_0'$ is a
relative weak mixing quotient of $\sigma$, in the sense of $2.9$ in
$\text{\rm [P1]}$.
\endproclaim
\noindent {\it Proof of Theorem} 4.1. Denote $\tilde{P}=P\overline{\otimes}
P$ and let $\sigma=\sigma_0\otimes \rho$,
$\tilde{\sigma}=\sigma_0\otimes \sigma_0 \otimes \rho$ be the
product actions of $\Gamma$ on $P\overline{\otimes} N$  and resp.
$\tilde{P} \overline{\otimes} N$.

Denote $M=P\overline{\otimes} N \rtimes \Gamma$,
$\tilde{M}=\tilde{P}\overline{\otimes} N \rtimes \Gamma$ and view
$M$ as a subalgebra of $\tilde{M}$ by identifying
$P\overline{\otimes} N$ with the subalgebra $(P\otimes 1)
\overline{\otimes} N$ of $P\overline{\otimes} P \overline{\otimes}
N=\tilde{P}\overline{\otimes} N$ and by identifying the canonical
unitaries $\{u_g\}_g$ in $M$, $\tilde{M}$ implementing $\sigma$ on
$P\overline{\otimes} N$ and $\tilde{\sigma}$ on
$\tilde{P}\overline{\otimes} N$. From now on we denote by $\tau$ the
canonical trace on the ambient algebra $\tilde{M}$ and on all its
subalgebras.

Since the s-malleable deformation $\alpha: \Bbb R \rightarrow
\text{\rm Aut}(\tilde{P},\tau)$ commutes with $\tilde{\sigma}$, it
extends to an action of $\Bbb R$ on $\tilde{M}$, still denoted
$\alpha$, equal to the identity on $N=1 \otimes N$ and on
$\{u_g\}_g$.

Note that if we denote $u'_g=w_gu_g$, then the cocyle relation for
$w_g$ is equivalent to the relation $u'_{g_1}u'_{g_2}=u'_{g_1g_2},
\forall g_1,g_2\in \Gamma$ in $M\subset \tilde{M}$. Also, denote by
$\sigma'$ the action of $\Gamma$ on $P\overline{\otimes} N$ given by
$\sigma'_g(x)=\text{\rm Ad}(u'_g)(x)=w_gu_gxu_g^*w_g^*
=w_g\sigma_g(x)w_g^*$. (N.B. If $P, N$ are commutative, then this is
equal to $\sigma_g(x)$.)

Note then that if we view $L^2 \tilde{M}$ as $L^2M
\overline{\otimes} L^2 P$ via the isomorphism $(x\otimes y) u_h
\mapsto (x u_h) \otimes y$, $x\in P\overline{\otimes} N$, $y\in P$,
$g\in \Gamma$, then the $\Gamma$-representation $\tilde{\pi}$ on
$L^2 \tilde{M}$ given by $\tilde{\pi}_g ( (x\otimes y) u_h) =
u'_g((x\otimes y) u_h){u'_g}^*$ corresponds to the
$\Gamma$-representation on $L^2 M \overline{\otimes} L^2 P$ given by
$\text{\rm Ad}(u'_g)(xu_h) \otimes \sigma_g (y)$, $\forall x\in
P\overline{\otimes}N, y\in P, g,h\in \Gamma$.

Thus, if $\sigma$ has stable spectral gap on $H$, then $\forall
\delta
> 0$, with $\delta < 2^{-5}$,
$\exists F\subset H$ finite and $\delta_0 > 0$ such that:  if
$u \in \Cal U(\tilde{M})$ satisfies $\|\tilde{\pi}_h(u)-u\|_2 \leq
\delta_0, \forall h\in F$, then $\|u-E_M(u)\|_2 \leq \delta$. Since
$\alpha_s(u'_h)$ is continuous in $s$ for all $h\in F$, it follows
that for sufficiently small $s>0$ we have
$\|\alpha_{-s/2}(u'_h)-u'_h\|_2 \leq \delta_0/2$, $\forall h\in F$.
Let $g$ be an arbitrary element in the group $H'$. Since the groups
$H,H'$ commute, $u'_{g}$ commutes with $u'_h, \forall h\in H$, in
particular $\forall h\in F$. Thus we get

$$
\|[\alpha_{s/2}(u'_{g}), u'_{h}]\|_2 =\|[u'_{g},
\alpha_{-s/2}(u'_h)]\|
$$
$$
\leq 2 \|\alpha_{-s/2}(u'_h)-u'_h\|_2 \leq \delta_0, \forall h\in F,
g\in H'.
$$

Since $\|[\alpha_{s/2}(u'_{g}), u'_{h}]\|_2
=\|\tilde{\pi}_h(\alpha_{s/2}(u'_{g}))-\alpha_{s/2}(u'_{g})\|_2$,
this implies that the unitaries $u=\alpha_{s/2}(u'_{g})\in
\tilde{M}, g\in H'$, satisfy the inequality
$\|\tilde{\pi}_h(u)-u\|_2 \leq \delta_0$. By the above conditions we
thus get

$$
\|\alpha_{s/2}(u'_{g})-E_M(\alpha_{s/2}(u'_{g}))\|_2 \leq \delta,
\forall g\in H',
$$
which by $(2.1.1)$ implies $\|\alpha_{s}(u'_{g})-u'_{g}\|_2 \leq
2\delta$, $\forall g\in H'$. Let $K=\overline{\text{\rm co}}^w
\{u'_{g}\alpha_{s}(u'_{g})^* \mid g\in H'\}$ and notice that $K$ is
a convex weakly compact subset, it is contained in the unit ball of
$\tilde{P}\overline{\otimes} N\subset \tilde{M}$ (because
$u'_{g}\alpha_{s}(u'_{g})^*=w_{h'}\alpha_{s}(w_{g})^*$) and for all
$\xi\in K$ and $g\in H'$ we have $u'_{g}\xi\alpha_{s}(u'_{g})^*\in
K$. Let $x\in K$ be the unique element of minimal norm
$\|\cdot\|_2$. Since $\|u'_g x \alpha_s(u'_g)^*\|_2=\|x\|_2, \forall
g\in H',$ by the uniqueness of $x$ it follows that
$u'_{g}x\alpha_{s}(u'_{g})^*=x$, $\forall g\in H'$. Thus $x$
intertwines the representations $g\mapsto u'_{g}$, $g \mapsto
\alpha_{s}(u'_{g})$. It follows that the partial isometry $v\in
\tilde{P}\overline{\otimes} N$ in the polar decomposition of $x$ is
non-zero and still intertwines the representations, i.e.
$u'_{g}v=v\alpha_{s}(u'_{g})$, or equivalently

$$
w_{g}\tilde{\sigma}_{g}(v)=v\alpha_{s}(w_{g}), \forall g\in H' \tag
4.1.1
$$

Moreover, since $\|u'_{g}\alpha_{s}(u'_{g})^*-1\|_2 =
\|u'_{g}-\alpha_{s}(u'_{g})\|_2\leq 2\delta$ we have $\|\xi-1\|_2
\leq 2 \delta, \forall \xi\in K$, thus $\|x-1\|_2 \leq 2 \delta$,
which by [C1] gives $\|v-1\|_2 \leq 4(2\delta)^{1/2}$.

By using the symmetry $\beta$, viewed as an automorphism of
$\tilde{P}\overline{\otimes} N$ (acting as the identity on $N$), the
same argument as in the proof of Lemma 4.6 in [P1] shows that
starting from $(4.1.1)$ applied to $s=2^{-n}$, for some large integer
$n$, one can obtain a partial isometry $v_1\in
\tilde{P}\overline{\otimes} N$ such that
$w_{g}\tilde{\sigma}_{g}(v_1)=v_1\alpha_1(w_{g}), \forall g\in H'$,
and $\|v_1\|_2 = \|v\|_2$. We repeat that argument below, for
completeness.

It is clearly sufficient to show that whenever we have $(4.1)$ for
some $s=2^{-n}$ and a partial isometry $v
\in\tilde{P}\overline{\otimes} N$, then there exists a partial
isometry $v'\in \tilde{P} \overline{\otimes} N$ satisfying
$\|v'\|_2=\|v\|_2$ and $w_g \tilde{\sigma}_g(v') = v'
\alpha_{2s}(w_g)$, $\forall g\in H'$. Indeed, because then the
statement follows by repeating the argument $n$ times.

Applying $\beta$ to $(4.1)$ and using the fact that $\beta$ commutes with
$\tilde{\sigma}$, $\beta(x)=x, \forall x\in P \overline{\otimes}
N\subset \tilde{P}\overline{\otimes} N$ and $\beta \alpha_s =
\alpha_{-s} \beta$ as automorphisms on $\tilde{P}\overline{\otimes}
N$, we get $\beta(w_h)=w_h$ and
$$
w_g \tilde{\sigma}_g(\beta(v))=\beta(v) \alpha_{-s}(w_g), \forall
g\in H'. \tag 4.1.2
$$
Since $(3.1)$ can be read as $v^* w_g =
\alpha_t(w_g)\tilde{\sigma}_g(v^*)$, from $(4.1.1)$ and $(4.1.2)$ we get
the identity

$$
v^* \beta(v) \alpha_{-s}(w_g)=v^* w_g \tilde{\sigma}_g(\beta(v))
$$
$$
=\alpha_s(w_g)\tilde{\sigma}_g(v^*)\tilde{\sigma}_g(\beta(v))
=\alpha_s(w_g)\tilde{\sigma}_g(v^*\beta(v)),
$$
for all $g\in H'$. By applying $\alpha_s$ on both sides of this
equality, if we denote $v'=\alpha_s(\beta(v)^*v)$, then we further
get
$$
{v'}^* w_g = \alpha_{2s} (w_g) \tilde{\sigma}_g({v'}^*), \forall
g\in H',
$$
showing that $w_g \tilde{\sigma}_g(v') = v' \alpha_{2s}(w_g)$,
$\forall g\in H'$, as desired. On the other hand, the intertwining
relation $(4.1.1)$ implies that $vv^*$ is in the fixed point algebra
$B$ of the action Ad$w_h \circ \tilde{\sigma}_g=\text{\rm Ad}(u'_g)$
of $H'$ on $\tilde{P}\overline{\otimes} N$. Since
$\tilde{\sigma}_{|H'}$ is weak mixing on $(1 \otimes P) \otimes 1
\subset \tilde{P}\overline{\otimes} N$ (because it coincides with
$\sigma_0$ on $1_P \otimes P \otimes 1_N \simeq P$) and because
Ad$w_h$ acts as the identity on $(1 \otimes P) \otimes 1$ and leaves
$(P\otimes 1) \overline{\otimes} N$ globally invariant, it follows
that $B$ is contained in $(P \otimes 1) \overline{\otimes} N$. Thus
$\beta$ acts as the identity on it (because it acts as the identity
on both $P\otimes 1$ and $1\otimes N$). In particular
$\beta(vv^*)=vv^*$, showing that the right support of $\beta(v^*)$
equals the left support of $v$. Thus, $\beta(v^*)v$ is a partial
isometry having the same right support as $v$, implying that $v'$ is
a partial isometry with $\|v'\|_2 = \|v\|_2$.

Altogether, this argument shows that $\forall \varepsilon_0 > 0$,
$\exists v_1 \in \tilde{P}\overline{\otimes} N$ partial isometry
satisfying $w_{g}\tilde{\sigma}_{g}(v_1)=v_1\alpha_1(w_{g}), \forall
g\in H'$, and $\|v_1\|_2\geq 1-\varepsilon_0/2$. Extending $v_1$ to
a unitary $u_1$ in $\tilde{P}\overline{\otimes} N$ it follows that
$\| w_{g}\tilde{\sigma}_{g}(u_1)-u_1\alpha_1(w_{g})\|_2 \leq
\varepsilon_0$, $\forall g\in H'$. By 2.12.2$^\circ$ in [P1] it
follows that the cocycles $w_g$ and $\alpha_1(w_g), g\in H'$, are
equivalent. Since $H'\curvearrowright X$ is assumed weak mixing, we
can apply Theorem 3.2 in [P1] to deduce that there exists $u\in
\Cal U(P \overline{\otimes} N)$ such that $w'_g=u^*w_g\sigma_g(u)$,
$g\in H',$ takes values into $\Cal U(N)$. By the weak mixing of $H'
\curvearrowright X$ and Lemma 3.6 in [P1], $w'_g$ takes values
into $\Cal U(N)$ for any $g$ in the w-normalizer of $H'$, in
particular on all $HH'$. The part of the statement concerning
wq-normalizer follows by applying again Lemma 3.6 in [P1], while
the part concerning actions $\sigma_0$ that extend to s-malleable
actions $\Gamma \curvearrowright^{\sigma_0'} P'$ such that that
$\sigma_0'$ is weak mixing relative to $\sigma_0$ follows from
Lemma 2.11 in [P1].

This ends the proof of Theorem 4.1. \hfill $\square$

\vskip .05in \noindent {\it Proof of Theorem} 1.1.
Let $(P,\tau) = (L^\infty
X, \int \text{\rm d}\mu)$, $(N,\tau_N)$ be a finite von Neumann
algebra such that $\Cal V$ is a closed subgroup of $\Cal U(N)$ and
$\rho=id$. If $w: X \times \Gamma \rightarrow \Cal V \subset \Cal
U(N)$ is a measurable (right) cocycle for $\Gamma \curvearrowright
X$, as defined for instance in 2.1 of [P1], then we view it as an
algebra (left) cocycle $w:\Gamma \rightarrow \Cal V^X \subset \Cal
U(L^\infty (X, N))=\Cal U(L^\infty X \overline{\otimes} N) = \Cal
U(P\overline{\otimes} N)$ for the action $\Gamma
\curvearrowright^\sigma  P\overline{\otimes} N$. The result follows
then from 4.1 and 3.5 in [P1]. \hfill $\square$ \vskip .05in
\noindent {\it Proof of Coroallry} 1.2. This is a trivial consequence
of Theorem 1.1 and Lemma
3.3. \vskip .05in \noindent {\it Proof of Corollary} 1.3. This
is an immediate consequence of
Theorem 1.1, Corollary 1.2 and 5.7-5.9 in [P1]. \vskip
.05in \noindent {\it Proof of Corollary} 1.4. This follows now readily from
Corollary 1.2 and 2.7 in [P2]. \hfill $\square$

\vskip .05in

We end this section by mentioning a non-commutative analogue of
Corollary 1.2 which follows from Theorem 4.1:

\proclaim{4.2. Corollary} Let $\Gamma$ be a countable group having
infinite commuting subgroups $H,H'$ with $H$ non-amenable. Let
$\Gamma \curvearrowright^\rho (N,\tau)$ be an arbitrary action. Let
$\Gamma \curvearrowright^{\sigma_0} (P,\tau)$ be an action whose
restriction to $HH'$ can be extended to an action $HH'
\curvearrowright^{\sigma'_0} (P',\tau)$ which is weak mixing
relative to ${\sigma_0}_{|HH'}$ and is one of the following: \vskip
.05in $1^\circ.$ A generalized non-commutative Bernoulli action $HH'
\curvearrowright (B_0,\tau_0)^I$, with base $(B_0,\tau_0)$ a finite
amenable von Neumann algebra and with the actions of $H,H'$ on the
countable set $I$ satisfying $|H' i|=\infty$ and $\{g \in H \mid
gi=i\}$ amenable, $\forall i\in I$.

$2^\circ.$ A Bogoliubov action associated to a unitary
representation of $HH'$ which has stable spectral gap on $H$ and no
finite dimensional $H'$-invariant subspaces. \vskip .05in If either
$H'$ is w-normal in $\Gamma$,  or if $H'$ is merely wq-normal in
$\Gamma$ but $\sigma_0$ is mixing, then any cocycle for
$\sigma_0 \otimes \rho$ is equivalent to a cocycle with values in
$N$.
\endproclaim

\heading 5. Proof of vNE Strong Rigidity  \endheading

To prove Theorem 1.5 we need two technical results about Bernoulli
actions, which are of independent interest. In fact,
these results hold
true for Bernoulli actions with arbitrary finite von Neumann
algebras as base, the proof being exactly the same as in the
commutative case.

Thus, we will denote by $(B_0,\tau_0)$ an amenable finite von
Neumann algebra and by $(B,\tau)$ a von Neumann subalgebra of
$(B_0,\tau_0)^\Lambda$ which is invariant to the Bernoulli action
$\Lambda \curvearrowright (B_0,\tau_0)^\Lambda$. Let
$\tilde{B}=B\overline{\otimes} B$, $M=B \rtimes \Lambda$,
$\tilde{M}=\tilde{B} \rtimes \Lambda$, where $\Lambda
\curvearrowright \tilde{B}$ is the double of the action $\Lambda
\curvearrowright B$. We view $M$ as a subalgebra of $\tilde{M}$ in
the obvious way, by identifying $B=B \otimes 1 \subset \tilde{B}$
and by viewing the canonical unitaries $\{v_h\mid h\in
\Lambda\}\subset \tilde{M}$ as also implementing $\Lambda
\curvearrowright B \otimes 1=B$.

\proclaim{5.1. Lemma} If $Q\subset M$ is a von Neumann subalgebra
with no amenable direct summand, then $Q'\cap \tilde{M}^\omega
\subset M^\omega$. Equivalently, $\forall \delta > 0$, $\exists F
\subset \Cal U(Q)$ finite and $\delta_0>0$ such that if $x\in
(\tilde{M})_1$ satisfies $\|[u,x]\|_2 \leq \delta_0$, $\forall u\in
F$, then $\|E_M(x)-x\|_2 \leq \delta$.
\endproclaim
{\it Proof}. By commuting squares of algebras, it is clearly
sufficient to prove the case $\Lambda \curvearrowright
(B,\tau)=(B_0,\tau_0)^\Lambda$. Let $\zeta_0=1, \zeta_1, ...$ be an
orthonormal basis of $L^2B_0$. Denote by $I$ the set of
multi-indices $(n_g)_g $ with $n_g \geq 0$, all but finitely many
equal to $0$. Note that $\Lambda$ acts on $I$ by left translation.
For each $i=(n_g)_g$ let $\eta_i=(\zeta_{n_g})_g$. Then
$\{\eta_i\}_i$ is an orthonormal basis of $L^2B,$ and $\Lambda
\curvearrowright B$ implements a representation $\Lambda
\curvearrowright L^2B$ which on $\xi_i$ is given by $g \xi_i=\xi_{g
i}$.

For each $i\in I_0=I \setminus \{0\}$, let $K_i=\{g \in \Lambda \mid
g i =i\}$ be the stabilizer of $i$ in $\Lambda$ and note that
$\Lambda i$ with the left translation by $\Lambda$ on it, is the same as
$\Lambda/K_i$. Denote $P_i = B \rtimes K_i\subset M$ and note that
since $K_i$ is a finite group $P_i$ is amenable. Let us show that
$L^2 ( \text{\rm sp} M (1 \otimes \xi_i) M, \tau)\simeq L^2 (\langle
M, e_{P_i} \rangle, Tr)$, as Hilbert $M$-bimodules. To this end, we
show that $x(1\otimes \xi_i)y\mapsto xe_{P_i}y$, $x,y\in M$, extends
to a well defined isomorphism between the two given Hilbert spaces.
It is in fact sufficient to show that
$$
\langle  x' (1\otimes \xi_i) y', x(1\otimes \xi_i) y \rangle_\tau
=\langle x' e_{P_i} y', xe_{P_i} y \rangle_{Tr},
$$
or equivalently
$$
Tr( y^*e_{P_i}x^*x'e_{P_i}y') = \tau( y^* (1 \otimes \xi_i^*) x^* x'
(1 \otimes \xi_i)y'), \tag a
$$
for all  $x,x', y,y'\in M$. Proving this identity for $x=v_g a$,
$x'=v_{g'}a'$, $y=v_h, y'=v_{h'}$, with $a,a'\in B = B \otimes 1$
and $g,g',h,h'\in \Lambda$ is clearly enough. The left side of $(a)$
is equal to
$$
\delta_{g^{-1}g',K_i}\tau(v_h^* a^* v_{g^{-1}g'} a'v_{h'}), \tag b
$$
where $\delta_{g^{-1}g',K_i}$ equals $0$ if $g^{-1}g'\not\in K_i$
and equals $1$ if $g^{-1}g'\in K_i$. On the other hand, the right
side of $(a)$ equals

$$
\tau(v_h^* a^* v_g^* (1 \otimes \xi_{-gi})(1\otimes \xi_{g' i})
v_{g'} a' v_{h'}) =\delta_{gi,g'i} \tau(v_h^* a^* v_{g^{-1}g'} a'
v_{h'}). \tag c
$$
Since $gi=g'i$ if and only if $g^{-1}g'\in K_i$, it follows that
$(b)=(c)$, showing that $(a)$ is indeed an identity.

We have thus shown that $L^2\tilde{M}\ominus L^2 M \simeq
\bigoplus_{i\in I_0} L^2 (\langle M, e_{P_i} \rangle, Tr)$. But since
$P_i$ are amenable, we have a weak containment of Hilbert
$M$-bimodules $L^2 (\langle M, e_{P_i} \rangle, Tr)$ $\prec L^2 M
\overline{\otimes} L^2 M$, $\forall i\in I_0$. Thus, we also have
such containment as Hilbert $Q$-bimodules.

On the other hand, if $Q'\cap \tilde{M}^\omega \not\subset
M^\omega$, then there exists a bounded sequence $(x_n)_n \in
\tilde{M}^\omega$ such that $E_M(x_n)=0$, $\|x_n\|_2=1$, $\forall
n$, and $\lim_n \|x_n y - yx_n \|_2=0$, $\forall y\in Q$. But this
implies $L^2 Q \prec L^2 \tilde{M}\ominus L^2 M$ as $Q$-bimodules.
From the above, this implies $L^2 Q \prec \bigoplus_{i\in I_0} L^2
(\langle M, e_{P_i} \rangle, Tr)\prec (L^2 M \overline{\otimes} L^2
M)^{I_0}$ as Hilbert $Q$-bimodules, which in turn shows that $Q$ has
a non-trivial amenable direct summand by Connes' Theorem (see the proof
of Lemma 2 in [P7]). \hfill $\square$

\vskip .05in

In the next lemma, the {\it w-normalizer} of a von Neumann
subalgebra $P_0\subset M$ is the smallest von Neumann subalgebra
$P\subset M$ that contains $P_0$ and has the property: if $uPu^*
\cap P$ is diffuse for some $u \in \Cal U(M)$, then $u\in P$.

\proclaim{5.2. Lemma} Assume that $\Lambda \curvearrowright
(B_0,\tau_0)^\Lambda$ is weak mixing relative to $\Lambda
\curvearrowright B$. Let $Q\subset pMp$ be a von Neumann subalgebra
with no amenable direct summand and with commutant $Q_0=Q'\cap pMp$
having no corner embeddable into $B$ inside $M$ $($e.g. if $B$ is
abelian, one can require $Q_0$ to be $\text{\rm II}_1$; in general
one can require $Q_0$ to have no amenable direct summand$)$. Then
there exists a non-zero partial isometry $v_0 \in M$ such that
$v_0^*v_0 \in Q_0'\cap pMp$ and $v_0Q_0v_0^* \subset L\Lambda$.
Moreover, if $\Lambda$ is ICC, then there exists a unitary element
$u\in \Cal U(M)$ such that $uQ_0u^*\subset L\Lambda$ and if $P$
denotes the w-normalizer algebra of $Q\vee Q_0$ in $pMp$, then
$uPu^*\subset L\Lambda$.
\endproclaim
{\it Proof}. It is clearly sufficient to prove the statement in case
$p=1$ (by taking appropriate amplifications of $Q\subset pMp$). We
may also clearly assume $\Lambda \curvearrowright
(B,\tau)=(B_0,\tau_0)^\Lambda$, by the relative weak mixing
condition (cf. [P4]). Moreover, we may assume the Bernoulli action
$\Lambda \curvearrowright B$ is s-malleable, i.e. $B_0=L^\infty
\Bbb T$ in the abelian case and $B_0=R$ in general. Indeed, because
any other abelian (resp. amenable) algebra $B_0$ can be embedded
into $L^\infty \Bbb T$ (resp. $R$) and $\Lambda \curvearrowright
(L^\infty \Bbb T)^\Lambda=L^\infty (\Bbb T^\Lambda)$ (resp. $\Lambda
\curvearrowright R^\Lambda$) is weak mixing relative to $\Lambda
\curvearrowright B_0^\Lambda$.

Let $\alpha: \Bbb R \rightarrow \text{\rm Aut}(\tilde{B})$, $\beta
\in \text{\rm Aut}(\tilde{B})$, $\beta^2=id$, give the s-malleable
path for the Bernoulli action $\Lambda \curvearrowright B$. Since
$\alpha, \beta$ commute with the double action $\Lambda
\curvearrowright \tilde{B}$, it follows that $\alpha$ (resp.
$\beta$) extends to an action, that we still denote by $\alpha$
(resp. $\beta$), of $\Bbb R$ (resp. $\Bbb Z/2\Bbb Z$) on
$\tilde{M}$.

We first prove that there exists a non-zero partial isometry $w\in
\tilde{M}$ such that $w^*w\in Q_0'\cap M$, $ww^* \in
\alpha_1(Q_0'\cap M)$, $wy=\alpha_1(y)w, \forall y\in Q_0$.

Fix $\varepsilon > 0$. By Lemma 5.1, there exists a finite set $F\subset
\Cal U(Q)$ and $\delta_0 > 0$ such that if $x\in (\tilde{M})_1$
satisfies $\|[u,x]\|_2 \leq \delta_0, \forall u\in F$, then
$\|E_M(x)-x\|_2 \leq \varepsilon/2$.

Since $\alpha_s(Q)$ commutes with $\alpha_s(Q_0)$ and $\alpha_s(u)$
is a continuous path, $\forall u\in F$, it follows that there exists
$n$ such that $s=2^{-n}$ satisfies
$$\|[u, \alpha_{s/2}(x)]\|_2=\|[\alpha_{-s/2}(u), x]\|_2
$$
$$
\leq 2 \|\alpha_{-s/2}(u)-u\|_2 \leq \delta_0, \forall x\in
(Q_0)_1, \forall u\in F.
$$
Thus $\|E_M(\alpha_{s/2}(x))-\alpha_{s/2}(x)\|_2 \leq \varepsilon/2$,
$\forall x\in (Q_0)_1$, in particular for all $x=v \in \Cal U(Q_0)$.
By $(2.1.1)$ and the choice of $\delta$, it follows that
$\|\alpha_{s}(v)-v\|_2 \leq \varepsilon$, $\forall v\in \Cal
U(Q_0)$. As in the proof of Theorem 1.1 in Section 4, this implies there
exists a partial isometry $V\in \tilde{M}$ such that $V v =
\alpha_s(v)V$, $\forall v\in \Cal U(Q_0)$ and $\|v-1\|_2 \leq 4
\varepsilon^{1/2}$. In particular, $V^*V\in Q_0'\cap \tilde{M}$,
$VV^* \in \alpha_s(Q_0'\cap \tilde{M})$ and $V \neq 0$ if
$\varepsilon < 1/16$. Since $Q_0$ has no corner that can be embedded
into $B$ inside $M$, by Theorem 3.2 in [P4] we have $Q_0'\cap
\tilde{M} = Q_0'\cap M$. But then exactly the same argument as in
the proof of 1.1 in Section 4 gives a partial isometry $V_1\in
\tilde{M}$ such that $\|V_1\|_2=\|V\|_2\neq 0$ and $V_1 v =
\alpha_1(v)V_1$, $\forall v\in \Cal U(Q_0)$.

By Steps 4 and 5 on page 395 in [P4], it then follows that there
exists a non-zero partial isometry $v_0 \in M$ such that $v_0^*v_0
\in Q_0'\cap M$ and $v_0Q_0v_0^* \subset L\Lambda$.

Assume now that $\Lambda$ is ICC, equivalently $L\Lambda$ is a
factor. As in the proof of 4.4 in [P4], to show that we can
actually get a unitary element $v_0$ satisfying $v_0Q_0v_0^* \subset
L\Lambda$, we use a maximality argument. Thus, we consider the set
$\Cal W$ of all families $(\{p_i\}_i, u)$ where $\{p_i\}_i$ are
partitions of 1 with projections in $Q_0'\cap M$, $u\in M$ is a
partial isometry with $u^*u=\Sigma_i p_i$ and $u(\Sigma_i Q_0p_i)u^*
\subset L\Lambda$. We endow $\Cal W$ with the order given by
$(\{p_i\}_i, u) \leq (\{p'_j\}_j, u')$ if $\{p_i\}_i \subset
\{p'_j\}_j$, $u=u'(\Sigma_i p_i)$. $(\Cal W, \leq)$ is clearly
inductively ordered.

Let $(\{p_i\}_i, u)$ be a maximal element. If $u$ is a unitary
element, then we are done. If not, then denote $q'= 1 -\Sigma_i p_i
\in Q_0'\cap M$ and take $q\in Q_0$ such that $\tau(qq')=1/n$ for
some integer $n \geq 1$. Denote $Q_1=M_{n \times n}(qQ_0qq')$
regarded as a von Neumann subalgebra of $M$, with the same unit as
$M$. Then the relative commutant of $Q_1$ in $M$ has no amenable
direct summand, so by the first part there exists a non-zero partial
isometry $w\in M$ such that $w^*w\in Q_1'\cap M$ and $wQ_1w^*
\subset L\Lambda$. Since $qq'\in Q_1$ has scalar central trace in
$Q_1$, it follows that there exists a non-zero projection in $w^*w
Q_1 w^*w$ majorized by $qq'$ in $Q_1$.

It follows that there exists a non-zero projection $q_0 \in
qq'Q_1qq' = qQ_0qq'$ and a partial isometry $w_0 \in M$ such that
$w_0^*w_0 = q_0$ and $w_0(qQ_0qq')w_0^* \subset L\Lambda$. Moreover,
by using the fact that $Q_0$ is diffuse, we may shrink $q_0$ if
necessary so that it is of the form $q_0 = q_1 q'\neq 0$ with $q_1
\in \Cal P(Q_0)$ of central trace equal to $m^{-1} z$ for some $z\in
\Cal Z(Q_0)$ and $m$ an integer. But then $w_0$ trivially extends to
a partial isometry $w_1 \in M$ with $w_1^*w_1 = q'z\in Q_0'\cap M$
and $w_1 Q_0w_1^* \subset L\Lambda$. Moreover, since $L\Lambda$ is a
factor, we can multiply $w_1$ from the left with a unitary element
in $L\Lambda$ so that $w_1w_1^*$ is perpendicular to $uu^*$. But
then $(\{p_i\}_i \cup \{q'z\}, u_1)$, where $u_1=u+w_1$, is clearly
in $\Cal W$ and is (strictly) larger than the maximal element
$(\{p_i\}_i, u)$, a contradiction.

We have thus shown that there exists a unitary element $u\in \Cal
U(M)$ such that $uQ_0u^*\subset L\Lambda$. But then by 3.1 in [P4]
it follows that $uQu^*\subset L\Lambda$ as well, and in fact all the
w-normalizer of $Q \vee Q_0$ is conjugated by $u$ into $L\Lambda$.
Thus, $uPu^*\subset L\Lambda$. \hfill $\square$

\vskip .05in \noindent {\it Proof of Theorem} 1.5. Let $H\subset \Gamma$ be
a non-amenable group with centralizer $H'=\{g\in \Gamma \mid gh=hg,
\forall h\in H\}$ non-virtually abelian and wq-normal in $\Gamma$.
With the above notations, we can take $B_0=L^\infty \Bbb T$. Let
$Q=\theta(L H) \subset pMp$ and $Q_0=\theta (LH)'\cap pMp$. By
hypothesis, $Q$ has no amenable direct summand and $Q_0$ is type
II$_1$. Thus, by Lemma 5.2 it follows that there exists $u\in \Cal
U(M)$ such that $uQ_0u^*\subset L\Lambda$. Moreover, since
$\theta(L\Gamma)$ is contained in the w-normalizer algebra $P$ of
$Q_0$, it follows that $u\theta(L\Gamma)u^*\subset L\Lambda$. From
this point on, the results in [P5] apply to conclude the proof.
\hfill $\square$

\vskip .05in \noindent {\it Proof of Theorem} 1.6. We may assume $L^\infty X
\rtimes \Gamma\subset L^\infty Y \rtimes \Lambda=M$, $L^\infty
X=L^\infty Y=A$ and for each $g\in \Gamma$ there exists a partition of
1 with projections $\{p^g_h\}_{h\in \Lambda}$ such that
$u_g=\Sigma_h p^g_hv_h$ give the canonical unitaries implementing
$\Gamma \curvearrowright A$. Thus, $Q=LH$ has no amenable
direct summand, $Q_0=L(H')$ is type II$_1$ and $L\Gamma$ is
contained in the w-normalizer algebra of $Q_0$.

By Lemma 5.2 it follows that there exists a unitary element $u$ in
$M$ such that $uL\Gamma u^* \subset L\Lambda$. Since $\Lambda
\curvearrowright A$ is Bernoulli, by Lemma 4.5 in [P5] it follows
that $\Gamma \curvearrowright A$ is mixing, thus Theorem 5.2 in
[P5]  applies to conclude that there exists $u\in \Cal U(M)$ such
that $u(\{u_h\}_h)u^* \subset \Bbb T \{v_h\}_h$ and $uAu^*=A$.
\hfill $\square$

\heading 6. Final remarks \endheading

\noindent {\bf 6.1. vNE versus OE: the Connes-Jones example}. Formal
definitions show that OE $\Rightarrow$ vNE, but are these notions of
equivalence really different, and if they are,
then how much different ? In other
words: If $\Gamma \curvearrowright X, \Lambda \curvearrowright Y$
are free ergodic m.p. actions, does $L^\infty X \rtimes \Gamma \simeq
L^\infty Y \rtimes \Lambda$ imply $\Cal R_\Gamma \simeq \Cal
R_\Lambda$? If $\theta$ denotes the isomorphism between the II$_1$
factors, this is same as asking whether there always exists $\rho\in
\text{\rm Aut}(L^\infty Y \rtimes \Lambda)$ such that
$\rho(\theta(L^\infty X))=L^\infty Y$.

Two sets of results give a positive answer to this question for
certain classes of group actions: On the one hand, if $\Gamma,
\Lambda$ are amenable, then by [OW] there does exist an
automorphism $\rho$ of $L^\infty Y \rtimes \Lambda\simeq R$ taking
$\theta(L^\infty X)$ onto $L^\infty Y$; in fact by [CFW] any two
{\it Cartan subalgebras} of $R$ are conjugated by an automorphism of
$R$. On the other hand, all vNE rigidity results in [P1,4,5,8],
[IPeP], [PV] are about showing that for any isomorphisms $\theta$
between certain group measure space factors $L^\infty X \rtimes
\Gamma, L^\infty Y \rtimes \Lambda$ (or even amplifications of such)
$\exists u\in \Cal U(L^\infty Y \rtimes \Lambda)$ such that
$\text{\rm Ad} u \circ \theta (L^\infty X) = L^\infty Y$. This is
unlike the amenable case though, where one can decompose $R$ in
uncountably many ways, $R=L^\infty X_i \rtimes_{\sigma_i} \Bbb Z$,
with $\Bbb Z \curvearrowright^{\sigma_i} X_i$ free ergodic actions
(which can even be taken cojugate to the same given $\Bbb
Z$-action), such that no inner automorphism of $R$ can take the
subalgebras $L^\infty X_i\subset R$ onto each other, for different
$i$'s ([FM]).

Nevertheless, the answer to ``vNE $\Rightarrow$ OE ?'' is negative
in general, as shown by Connes and Jones in
[CJ1] through the following example:
Let $\Gamma_0$ be any non-amenable group and
$\Gamma_1=\Sigma_n H_n$ an infinite direct sum of non-abelian groups.
Let $H_n \curvearrowright [0,1]$ be any free $H_n$-action preserving
the Lebesgue measure (e.g. a Bernoulli $H_n$-action) and let
$\Gamma_1=\Sigma_n H_n\curvearrowright [0,1]^{\Bbb N}$ be the
product of these actions. Finally, denote $X=([0,1]^{\Bbb
N})^{\Gamma_0}$ and let $\Gamma_0$ act on $X$ by (left) Bernoulli
shifts and $H_0$ act diagonally, identically on each copy of
$[0,1]^{\Bbb N}$. Since the $\Gamma_0, \Gamma_1$ actions commute
they implement an action of $\Gamma= \Gamma_0 \times \Gamma_1$ on
$(X,\mu)$, which is easily seen to be free.

Since $\Gamma_0$ is non-amenable and $\Gamma_0 \curvearrowright X$
is Bernoulli, $\Gamma \curvearrowright X$ is strongly ergodic (it
even has spectral gap), thus $\Cal R_\Gamma$ is strongly ergodic as
well. However, since any sequence of canonical unitaries $v_{h_n}$
with $h_n \in H_n$ is central for $M=L^\infty X \rtimes \Gamma$, by
the non-commutativity of the $H_n$'s it follows that $M'\cap
M^\omega$ is non-commutative, so by McDuff's theorem $M\simeq
M\overline{\otimes} R$. Thus $M$ can also be decomposed as
$M=(L^\infty X \rtimes \Gamma)\overline{\otimes}(L^\infty ([0,1])
\rtimes H)=L^\infty (X \times [0,1]) \rtimes (\Gamma \times H)$,
where $H\curvearrowright [0,1]$ is any free ergodic m.p. action of an
amenable group $H$. Such $H \curvearrowright L^\infty ([0,1])$ always has
non-trivial approximately invariant sequences, i.e. it is not
strongly ergodic. Thus $\Gamma \times H \curvearrowright
L^\infty (X\times [0,1])$ is not strongly ergodic either, so it
cannot be OE to $\Gamma \curvearrowright X$ although both actions
give the same II$_1$ factor, i.e. are vNE. Thus vNE $ \not\Rightarrow$
OE.

Note that by taking $H_n=H$ , $\forall n$, one gets the same group
$\Gamma\simeq \Gamma \times H$ having two actions, one strongly
ergodic the other not, both giving rise to the same II$_1$ factor.
Moreover, if $\Gamma_0$ is taken Kazhdan, or merely w-rigid, then
the action $\Gamma_0 \curvearrowright X$ satisfies the hypothesis of
5.2/5.3 in [P1], so it is cocycle superrigid. Since it is weakly
mixing, its extension to $\Gamma \curvearrowright X$ is also cocycle
superrigid. Similarly, if $\Gamma_0$ is taken as a product between a
non-amenable group and an infinite group, then $\Gamma
\curvearrowright X$ follows cocycle superrigid by Theorem 1.1. In
particular, in both cases $\Gamma \curvearrowright X$ is OE
superrigid, so by 5.7 in [P1] and
Corollary 1.3 $\mycal F(\Cal R_\Gamma)$ is
countable. If in addition $\Gamma_0$ and $H_n$ have no finite normal
subgroups, $\forall n$, then $\mycal F(\Cal R_\Gamma)=1$.

In other words, there exists a free ergodic cocycle superrigid
action $\Gamma \curvearrowright X$ which is strongly ergodic,
satisfies $\mycal F(\Cal R_\Gamma)=1$, but the associated II$_1$
factor $M=L^\infty X \rtimes \Gamma$ can also be realized as
$M=L^\infty Y \rtimes \Gamma'$ with $\Gamma' \curvearrowright Y$
a free ergodic but not strongly ergodic action with $\Cal R_{\Gamma'}
\simeq \Cal R_{\Gamma'} \times \Cal R_{hyp}$, $M \simeq
M\overline{\otimes} R$. In particular $\mycal F(M) =\mycal F(\Cal
R_{\Gamma'})=\Bbb R_+^*$. Moreover, one can take $\Gamma \simeq
\Gamma'$.

\vskip .1in

\noindent {\bf 6.2. On the transversality of mallebale actions}.
Although all existing examples of malleable actions are in fact
s-malleable, it would be interesting to give a proof of Theorem 1.1
that would only use (basic) malleability, even if this means
sacrificing some of the generality on the side of the target groups.
For instance, to prove Theorem 1.1 for cocycles of malleable actions with
abelian, compact or discrete groups as targets. But it seems to us
that any alternative argument would still need some sort of
``transversality'' property for an appropriate family
$\{\alpha_s\}_s$ of automorphisms commuting with the double action
$\Gamma \curvearrowright X \times X$ and relating $id$ to the flip,
requiring that if $\alpha_s(x)$ close to $L^\infty X \otimes 1$ for
some $x\in L^\infty X \otimes 1$, then $\alpha(x)$ is close to $x$.
Besides s-malleability, another sufficient condition for this to
happen is the following:

\vskip .05in \noindent $(6.2.1)$ There exists a Hilbert space $\Cal K$
containing $L^2(X \times X, \mu \times \mu)$, an orthonormal system
$\{\xi_n\}_n$ $\subset \Cal K$ satisfying $L^2 X \otimes 1  \subset
\overline{\Sigma_n \Bbb C \xi_n}$, and an extension of $\alpha_s$ to
a unitary element $\alpha'_s$ on $\Cal K$, such that $\langle \xi_n,
\alpha'_s(\xi_m)\rangle = \delta_{nm} c_n$,  with $c_n \in \Bbb R$,
$\forall n,m$. \vskip .05in

Indeed, it is easy to see that if an automorphism $\alpha_s$
satisfies $(6.2.1)$, then $\|\alpha_s^2 (x)-x\|_2 \leq 2 \sqrt{2}
\|\alpha_s(x)- E_{L^\infty X} (\alpha_s(x))\|_2$, $\forall x\in
L^\infty X \otimes 1$. In fact, in an initial version of this paper
we used property $(6.2.1)$ to derive the transversality $(2.1.1)$, and
proved that Bernoulli, Gaussian and Bogoliubov actions satisfy
$(2.1.1)$ by showing they satisfy $(6.2.1)$. It was Stefaan Vaes and the
referee who pointed out to us that in fact s-malleability trivially
implies the transversality condition $(2.1.1)$ (i.e. Lemma 2.1).

Nevertheless, condition $(6.2.1)$ seems interesting in its own right.
Related to it, note that if $\Gamma \curvearrowright X$ is so that
Aut$_\Gamma (X \times X)$ contains a finite group $K$ that has the
flip in it and for which there exists an extension of $K
\curvearrowright L^2 X \otimes L^2 X$ to a representation $K
\curvearrowright \Cal K$, with an orthonormal system
$\{\xi_n\}_n\subset \Cal K$ spanning $L^2 X \otimes 1$, such that
the Hilbert spaces $\Cal K_n=\text{\rm sp} \{k \xi_n \mid k \in K\}$
are mutually orthogonal and have dimensions majorized by some
constant $c=c(|K|)$ with the property that $\forall n$, $\exists k
\in K\setminus \{e\}$ with $\|k\xi_n - \xi_n\|_2 < 1$, then $\Gamma
\curvearrowright X$ would automatically satisfy a cocycle
superrigidity result, with no additional requirements on the group
$\Gamma$, or on the way it acts on $X$.

\vskip .1in \noindent {\bf 6.3. $\Cal C\Cal S$ and} $\Cal O\Cal
E\Cal S$ {\bf groups}. Related to Remark 6.7 in [P1], we
re-iterate here the following question: What is the class $\Cal
C\Cal S$ of groups $\Gamma$ for which the Bernoulli action $\Gamma
\curvearrowright \Bbb T^\Gamma$ is $\Cal U_{fin}$-cocycle superrigid
? (N.B. Any relative weak mixing quotient of $\Gamma
\curvearrowright \Bbb T^\Gamma$, for $\Gamma \in \Cal C\Cal S$, is
then automatically $\Cal U_{fin}$-CSR as well, by results in [P1].)
The class $\Cal C\Cal S$ cannot contain free products with
amalgamation $\Gamma=\Gamma_1 *_H \Gamma_2$, with $H$ a finite
subgroup of $\Gamma_i$, $H\neq \Gamma_i$, $i=1,2$ (see e.g. [P2]).
The class covered by Theorem 1.1 does not contain word hyperbolic
groups. Hyperbolic groups with Haagerup property are not covered by
5.2/5.3 in [P1] either, because they cannot have infinite
subgroups with the relative property (T).

The following question is equally interesting: What is the class of
groups $\Gamma$ for which any OE between a Bernoulli $\Gamma$-action
$\Gamma \curvearrowright (X_0,\mu_0)^\Gamma$ and an arbitrary
Bernoulli action $\Lambda \curvearrowright (Y_0,\nu_0)^\Lambda$
comes from a conjugacy. It is very possible that this class consists
of all non-amenable groups. It would be very interesting to decide
this question for the free groups. A related question is to
characterize the sub-class $\Cal O\Cal E\Cal S$ of groups $\Gamma$
for which the Bernoulli action $\Gamma \curvearrowright \Bbb
T^\Gamma$ is OE Superrigid. $\Cal O\Cal E\Cal S$ doesn't contain any
free product of infinite amenable groups, by [OW], [CFW].

\vskip .1in \noindent {\bf 6.4. Examples of prime factors}. Lemma
5.2 allows deriving new examples of prime II$_1$ factors, i.e.
factors $M$ that cannot be decomposed as  tensor products $M=Q
\overline{\otimes} Q_0$ with $Q,Q_0$ II$_1$ factors (see [O1], [O2],
[Pe] for other examples of such factors):

\proclaim{6.4.1. Theorem} Let $\Lambda$ be an arbitrary non-amenable
group and $\Lambda \curvearrowright Y$ a free relative weak mixing
quotient of a Bernoulli action. Then $L^\infty Y \rtimes \Lambda$ is
prime. More generally, if $B\subset R^\Lambda$ is a von Neumann
algebra invariant to the action $\Lambda\curvearrowright R^\Lambda$,
such that $\Lambda \curvearrowright B$ is free and $\Lambda
\curvearrowright R^\Lambda$ is weak mixing relative to $\Lambda
\curvearrowright B$, then $B \rtimes \Lambda$ is prime. In
particular $L^\infty \Bbb T^\Lambda \rtimes \Lambda$ and $R^\Lambda
\rtimes \Lambda$ are prime.
\endproclaim
\noindent {\it Proof}. Denote $M=L^\infty Y \rtimes
\Lambda$. Assume $M = Q \overline{\otimes} Q_0$. Since $M$ is
non$(\Gamma)$ (see e.g. [I1]), it follows that both $Q,Q_0$ are
non$(\Gamma)$, thus non-amenable. By the first part of Lemma 5.2,
there exists a non-zero $p\in Q_0'\cap M=Q$ and a unitary element
$u\in M$ such that $u(Q_0p)u^*\subset L\Lambda$. By 3.1 in [P4] it
follows that $up( Q \vee Q_0)pu^*\subset L\Lambda$. But the left
hand side is equal to $p'Mp'$, where $p'=upu^*$. This means
$p'L\Lambda p'=p'Mp'$, a contradiction. \hfill $\square$

\vskip .05in

We mention that a more careful handling of the proof of Lemmas 5.1,
5.2 allows us to prove that factors $B \rtimes \Lambda$ associated to
Bernoulli actions $\Lambda \curvearrowright (B, \tau)= (B_0,
\tau_0)^\Lambda$, with an arbitrary finite von Neumann algebra $B_0\neq
\Bbb C$ as base, are prime for any non-amenable $\Lambda$ (see [I2]
for related rigidity results on such factors).

Note that Lemmas 5.1, 5.2 show that if $\Lambda$ is an ICC group such that
$M=L\Lambda$ has the property:

\vskip .05in \noindent $(6.4.1)$ {\it If $Q\subset M$ has type}
II$_1$ {\it relative commutant $Q'\cap M$, then $Q$ is amenable},
\vskip .05in \noindent then given any free, relative weak mixing
quotient $\Lambda \curvearrowright Y$ of the Bernoulli action
$\Lambda \curvearrowright \Bbb T^\Lambda$, the II$_1$ factor
$M=L^\infty Y \rtimes \Lambda$ has property $(6.4.1)$ as well.
Indeed, because if $Q\subset M$ has no amenable direct summand and
$Q_0=Q'\cap M$ is of type II$_1$, then by the last part of Lemma 5.2 there
exists a unitary element $u\in M$ such that $u(Q \vee Q_0)u^*\subset
L\Lambda$, contradicting the property for $L\Lambda$. This result
should be compared with a result in [O2], showing that if
$\Lambda$ satisfies property AO and $H$ is an abelian group, then the
wreath product $H \wr \Lambda$ has the property AO as well. By
[O1] this implies $L(\Lambda \ltimes H) = L^\infty \hat{H} \rtimes
\Lambda$ is solid, thus prime.

\vskip .1in \noindent {\bf 6.5. On spectral gap rigidity}. The
results in Theorems 1.1
through 1.6 add to the plethora of rigidity phenomena involving
product groups that have been discovered in recent years in group
theory, OE ergodic theory, Borel equivalence relations  and von
Neumann algebras/II$_1$ factors ([MoSh], [HK], [OP], [Mo], [BSh],
etc). It would of course be interesting to find some common ground
(explanation) to these results. The idea behind our approach is very
much in the spirit of II$_1$ factor theory, but is otherwise rather
elementary. It grew out from an observation in [P6], where for the
first time spectral gap rigidity was used to prove a structural
rigidity result for II$_1$ factors. The starting point of all
deformation/spectral gap rigidity arguments we have used in this
paper and in [P6], [P7] is the following observation, which can be
viewed as a general ``spectral gap rigidity principle'':

\proclaim{6.5.1. Lemma} Let $\Cal U$ be a group of unitaries in a
$\text{\rm II}_1$ factor $\tilde{M}$ and $M, \tilde{P} \subset
\tilde{M}$ von Neumann subalgebras such that $\Cal U$ normalizes
$\tilde{P}$ and the commutant of $\Cal U$ in $\tilde{P}$, $Q_0=\Cal
U'\cap \tilde{P}$, is contained in $M$. Assume:

\vskip .05in \noindent $(6.5.1)$ The action $\text{\rm Ad} \Cal U$
on $\tilde{P}$ has spectral gap relative to $M$, i.e. for any
$\varepsilon > 0$, there exist
$F(\varepsilon) \subset \Cal U$ finite and
$\delta(\varepsilon)
> 0$ such that if $x\in (\tilde{P})_1$, $\|uxu^*-x\|_2 \leq \delta(\epsilon)$,
$\forall u\in F(\epsilon)$, then $\|E_{M}(x)-x\|_2 \leq \varepsilon$.
$($Note that this is equivalent to the condition $\Cal U'\cap
\tilde{P}^\omega\subset M^\omega)$.

\vskip .05in

Then any deformation of $id_{\tilde{M}}$ by automorphisms $\theta_n
\in \text{\rm Aut}(\tilde{M})$ satisfies:
$$
\lim_n (\sup\{ \|\theta_n(y)-E_M(\theta_n(y))\|_2 \mid y\in
(Q_0)_1\})=0. \tag 6.5.2
$$
In other words, the unit ball of $\theta_n(Q_0)$ tends to be
contained into the unit ball of $M$, as $n \rightarrow \infty$.
\endproclaim
\noindent {\it Proof}. Fix $\varepsilon > 0$ and let $F(\varepsilon)
\subset \Cal U$, $\delta(\varepsilon) > 0$, as given by $(6.5.1)$.
Let $n$ be large enough so that $\|\theta_n(u)-u\|_2 \leq \delta/2$,
$\forall u\in F$. If $x \in (\theta_n(Q_0))_1$ then $x$ commutes
with $\theta_n(F)$ and thus $\|uxu^*-x\|_2 \leq 2
\|u-\theta_n(u)\|_2 \leq \delta$. By $(6.5.1)$, this implies
$\|x-E_M(x)\|_2 \leq \varepsilon$. \hfill $\square$

\vskip .05in

In the proof of the cocycle superrigidity
result in Theorem1.1, Lemma 6.5.1 is used for
$\tilde{P}=L^\infty X\overline{\otimes} L^\infty X
\overline{\otimes}N$, $\tilde{M}=\tilde{P} \rtimes \Gamma$,
$M=(L^\infty X \otimes 1 \overline{\otimes}N) \rtimes \Gamma$ and
$\Cal U=\{u_h\mid h\in H\}$.

In the proof of the vNE and OE strong rigidity
results in Theorems 1.5 and 1.6, Lemma 6.5.1 is
used for $\tilde{M}=\tilde{P}=L^\infty Y \overline{\otimes} L^\infty
Y \rtimes \Lambda$, $M=L^\infty Y \rtimes \Lambda$ and $\Cal
U=\theta(\{u_h\mid h \in H\})$.

In the proof of Theorem 1 in [P7] it is used for
$\tilde{M}=\tilde{P}=L\Bbb F_n * L\Bbb F_n$, $M=L\Bbb F_n * \Bbb C$,
$\Cal U=\Cal U(Q)$.

In all these cases the deformation of $id_{\tilde{M}}$ is by
automorphisms of a malleable path $\alpha_s, s\in \Bbb R$.

The initial result in [P7], where a ``baby version'' of spectral
gap rigidity was used, states that if $\tilde{M}=Q\overline{\otimes}
R$ is a McDuff II$_1$ factor, with $Q$ non$(\Gamma)$, then any other
tensor product decomposition $\tilde{M}=N\overline{\otimes}P$ with
$N$ non$(\Gamma)$ and $P \simeq R$ is unitary conjugate to it, after
re-scaling. In this case one applies Lemma 6.5.1 for
$\tilde{P}=\tilde{M}$, $\Cal U=\Cal U(Q)$, $M=Q_0=R$. The trick then
is to take a deformation by inner automorphisms $\theta_n=\text{\rm
Ad}(v_n)$ with $v_n \in \Cal U(R_n)$ where $R_n \subset R$ is a
decreasing sequence of subfactors splitting off the $2^n$ by $2^n$
matrices in $R$, i.e. $R=R_n \otimes M_{2^n \times 2^n}(\Bbb C)$,
and satisfying $\bigcap_n R_n = \Bbb C1$. By Lemma 6.5.1 one then gets
$vQ_0v^* \approx Q_0$ (unit balls) uniformly in $v \in \Cal U(R_n)$,
for $n$ large, implying that $\Cal U(R_n) \underset{\sim} \to
\subset Q_0$, thus $R_n \underset{\sim} \to \subset Q_0$ (unit
balls), so by [OP] there exists $u\in \Cal U(\tilde{M})$ with the
required properties. Note that there is an alternative way to carry
out this argument, using the deformation by conditional expectations
$E_{R_n'\cap \tilde{M}}$, as explained in \S 5 of [P7].

\head  References\endhead

\item{[BSh]} U. Bader, Y. Shalom: {\it Factor and normal subgroup
theorems for lattices in products of groups}, preprint 2005,
to appear in Invent. Math.

\item{[CCJJV]} P.-A. Cherix, M. Cowling, P. Jolissaint,
P. Julg, A. Valette: ``Groups
with Haagerup property'', \newline Birkh$\ddot{\text{\rm a}}$user
Verlag, Basel Berlin Boston, 2000.

\item{[Ch]} E. Christensen: {\it Subalgebras of a finite algebra}, Math.
Ann. {\bf 243} (1979), 17-29.

\item{[C1]} A. Connes: {\it Classification of injective factors},
Ann. of Math., {\bf 104} (1976), 73-115.

\item{[C2]} A. Connes: {\it Sur la classification des facteurs
de type} II, C. R. Acad. Sci. Paris {\bf 281} (1975), 13-15.

\item{[C3]} A. Connes: {\it A type II$_1$ factor with countable
fundamental group}, J. Operator Theory {\bf 4} (1980), 151-153.

\item{[C4]} A. Connes: {\it Classification des facteurs},
Proc. Symp. Pure Math. {\bf 38}, Amer. Math. Soc. 1982, 43-109.

\item{[CFW]} A. Connes, J. Feldman, B. Weiss: {\it An amenable equivalence
relation is generated by a single transformation}, Erg. Theory Dyn.
Sys. {\bf 1} (1981), 431-450.

\item{[CJ1]} A. Connes, V.F.R. Jones: {\it A} II$_1$ {\it factor
with two non-conjugate Cartan subalgebras}, Bull. Amer. Math. Soc.
{\bf 6} (1982), 211-212.

\item{[CJ2]} A. Connes, V.F.R. Jones: {\it Property} (T)
{\it for von Neumann algebras}, Bull. London Math. Soc. {\bf 17}
(1985), 57-62.

\item{[CW]} A. Connes, B. Weiss: {\it Property} (T) {\it and
asymptotically invariant sequences}, Israel J. Math. {\bf 37}
(1980), 209-210.

\item{[Dy]} H. Dye: {\it On groups of measure preserving
transformations}, II, Amer. J. Math, {\bf 85} (1963), 551-576.

\item{[FM]} J. Feldman, C.C. Moore: {\it Ergodic equivalence
relations, cohomology, and von Neumann algebras I, II}, Trans. Amer.
Math. Soc. {\bf 234} (1977), 289-324, 325-359.

\item{[Fu1]} A. Furman: {\it Orbit equivalence
rigidity}, Ann. Math. {\bf 150} (1999), 1083-1108.

\item{[Fu2]} A. Furman: {\it On Popa's Cocycle Superrigidity
Theorem}, math.DS/0608364, \newline preprint 2006.

\item{[F]} H. Furstenberg: {\it
Ergodic behavior of diagonal measures and a theorem of Szemeredi on
arithmetic progressions}, J. d'Analyse Math. {\bf 31} (1977)
204-256.

\item{[G1]} D. Gaboriau: {\it Cout des r\'elations d'\'equivalence
et des groupes}, Invent. Math. {\bf 139} (2000), 41-98.

\item{[G2]} D. Gaboriau: {\it Invariants $\ell^2$ de r\'elations
d'\'equivalence et de groupes},  Publ. Math. I.H.\'E.S. {\bf 95}
(2002), 93-150.

\item{[GP]} D. Gaboriau, S. Popa: {\it An Uncountable Family of
Non Orbit Equivalent Actions of $\Bbb F_n$}, Journal of AMS
{\bf 18} (2005), 547-559.

\item{[Ha]} U. Haagerup: {\it An example of a non-nuclear}
$C^*$-{\it algebra which has the metric approximation property},
Invent. Math. {\bf 50} (1979), 279-293.

\item{[H]} G. Hjorth: {\it A converse to Dye's theorem}, Trans AMS {\bf 357}
(2004), 3083-3103.

\item{[HK]} G. Hjorth, A. Kechris: ``Rigidity theorems
for actions of product groups and
countable Borel equivalence relations'', Memoirs of the Amer. Math.
Soc., {\bf 177}, No. 833, 2005.

\item{[I1]} A. Ioana: {\it A relative version of Connes} $\chi(M)$
{\it invariant}, Ergod. Th. \& Dynam.
Sys. {\bf 27} (2007), 1199-1213 (math.OA/0411164).

\item{[I2]} A. Ioana: {\it Rigidity results for wreath product}
II$_1$ {\it factors}, math.OA/0606574.

\item{[I3]} A. Ioana: {\it Existence of uncountable families of
orbit inequivalent actions for groups containing $\Bbb F_2$},
preprint 2007.

\item{[IPeP]} A. Ioana, J. Peterson, S. Popa: {\it Amalgamated free
products of w-rigid factors and calculation of their symmetry
groups}, math.OA/0505589, to appear in Acta Math.

\item{[K]} D. Kazhdan: {\it Connection of the dual space of a
group with the structure of its closed subgroups}, Funct. Anal. and
its Appl. {\bf 1} (1967), 63-65.

\item{[Mc]} D. McDuff: {\it Central sequences and the hyperfinite
factor}, Proc. London Math. Soc. {\bf 21} (1970), 443-461.

\item{[Mo]} N. Monod: {\it Superrigidity for irreducible
lattices and geometric splitting},
Journal Amer. Math. Soc.,  {\bf 19} (2006) 781-814

\item{[MoSh]} N. Monod, Y. Shalom:
{\it Orbit equivalence rigidity and bounded cohomology}, Annals of Math.,
{\bf 164} (2006).

\item{[MvN1]} F. Murray, J. von Neumann:
{\it On rings of operators}, Ann. Math. {\bf 37} (1936), 116-229.

\item{[MvN2]} F. Murray, J. von Neumann: {\it Rings of operators
IV}, Ann. Math. {\bf 44} (1943), 716-808.

\item{[OW]} D. Ornstein, B. Weiss: {\it Ergodic theory of
amenable group actions I. The Rohlin Lemma} Bull. A.M.S. (1) {\bf 2}
(1980), 161-164.

\item{[O1]} N. Ozawa: {\it Solid von Neumann algebras},
Acta Math. {\bf 192} (2004), 111-117.

\item{[O2]} N. Ozawa: {\it A Kurosh type theorem for type} II$_1$
{\it factors}, Int. Math. Res. Notices, 2006,  math.OA/0401121

\item{[OP]} N. Ozawa, S. Popa: {\it Some prime factorization results
for type} II$_1$ {\it factors}, Invent Math. {\bf 156} (2004),
223-234.

\item{[Pe]} J. Peterson, {\it $L^2$-rigidity in von Neumann algebras},
math.OA/0605033.

\item{[P1]} S. Popa: {\it Cocycle and orbit equivalence
superrigidity for malleable actions of $w$-rigid groups},
Invent. Math. First on Line DOI 10.1007/s00222-007-0063-0
(math.GR/0512646).

\item{[P2]} S. Popa: {\it Some computations of $1$-cohomology groups
and construction of non orbit equivalent actions}, Journal of the
Inst. of Math. Jussieu {\bf 5} (2006), 309-332.

\item{[P3]} S. Popa: {\it Some rigidity results for
non-commutative Bernoulli shifts}, J. Fnal. Analysis {\bf 230}
(2006), 273-328.

\item{[P4]} S. Popa: {\it Strong Rigidity of} II$_1$ {\it Factors
Arising from Malleable Actions of $w$-Rigid Groups} I, Invent. Math.
{\bf 165} (2006), 369-408 (math.OA/0305306).

\item{[P5]} S. Popa: {\it Strong Rigidity of} II$_1$ {\it Factors
Arising from Malleable Actions of $w$-Rigid Groups} II, Invent.
Math. {\bf 165} (2006), 409-452 (math.OA/0407137).

\item{[P6]} S. Popa: {\it Deformation and rigidity in the study of} II$_1$
{\it factors}, Mini-Course at College de France, Nov. 2004.

\item{[P7]} S. Popa: {\it On Ozawa's property for free group
factors}, Math. Res. Notices.
Vol. {\bf 2007}, article ID rnm036, 10 pages,
doi:10.1093/imrn/rnm036 published on June 22, 2007 (math.OA/0607561).

\item{[P8]} S. Popa: {\it On a class of type II$_1$ factors with
Betti numbers invariants}, Ann. Math. {\bf 163} (2006), 809-889
(math.OA/0209310).

\item{[P9]} S. Popa: {\it Deformation and rigidity for group actions
and von Neumann algebras}, in ``Proceedings of the International
Congress of Mathematicians'' (Madrid 2006), Volume I, EMS Publishing House,
Zurich 2006/2007, pp. 445-479.

\item{[PSa]} S. Popa, R. Sasyk:
{\it On the cohomology of Bernoulli actions}, Erg. Theory Dyn. Sys.
{\bf 26} (2006), 1-11 (math.OA/0310211).

\item{[PV]} S. Popa, S. Vaes: {\it Strong rigidity of generalized Bernoulli
actions and computations of their symmetry groups}, math.OA/0605456,
to appear in Adv. in Math.

\item{[PoSt]} R. Powers, E. St\o rmer: {\it Free states of the
canonical anticommutation relations}, Comm. Math. Phys. {\bf 16}
(1970), 1-33.

\item{[Sc]} K. Schmidt: {\it Asymptotically invariant sequences
and an action of $SL(2, \Bbb Z)$ on the $2$-sphere}, Israel. J.
Math. {\bf 37} (1980), 193-208.

\item{[Sh]} Y. Shalom: {\it Measurable group theory}, In
``European Congress of Mathematics'' (Stockholm 2004), European Math
Soc, Zurich 2005, 391-424.

\item{[Si]} I.M. Singer: {\it Automorphisms of finite factors},
Amer. J. Math. {\bf 77} (1955), 117-133.

\item{[V]} S. Vaes: {\it Rigidity results for Bernoulli
actions and their von Neumann algebras} (after Sorin Popa)
S\'eminaire Bourbaki, exposé 961. Ast\'erisque (to appear).

\item{[Z1]} R. Zimmer: {\it Strong rigidity for ergodic actions of
seimisimple Lie groups}, Ann. of Math. {\bf 112} (1980), 511-529.

\item{[Z2]} R. Zimmer: ``Ergodic Theory and Semisimple Groups'',
Birkhauser, Boston, 1984.

\item{[Z3]} R. Zimmer: {\it Extensions of ergodic group actions},
Illinois J. Math. {\bf 20} (1976), 373-409.

\enddocument